\documentclass[11pt]{article}
\usepackage{amssymb}
\usepackage{amsmath}
\usepackage{latexsym}
\usepackage{amsfonts}
\usepackage{amsthm}
\usepackage[all]{xy}
\usepackage{graphicx}

\theoremstyle{plain}
\newtheorem{thm}{Theorem}
\newtheorem*{thm*}{Theorem}

\newtheorem{defn}{Definition}
\newtheorem{lem}{Lemma}
\newtheorem{prop}{Proposition}
\newtheorem{cor}{Corollary}

\theoremstyle{definition}
\newtheorem{ex}{Example}

\newtheorem*{dem}{Proof}

\theoremstyle{remark}
\newtheorem*{rmk}{\underline{Remark}}

\begin{document}
\title{Cohomological invariants of finite Coxeter groups}
\author {J\'er\^ome DUCOAT\\\\\\
e-mail : jerome.ducoat@ujf-grenoble.fr}
\maketitle

\begin{abstract}In this paper, we generalize Serre's splitting principle for cohomological invariants of the symmetric group to finite Coxeter groups, provided that the ground field has characteristic zero. We then determine all the cohomological invariants of the Weyl groups of classical type with coefficients in $\mathbb Z/2\mathbb Z$.
\end{abstract}  

\section*{Introduction}
Let $k_0$ be a field. We denote by $\Gamma_{k_0}$ the absolute Galois group of $k_0$.
Let $G$ be a smooth linear algebraic group scheme over $k_0$ and let $C$ be a finite $\Gamma_{k_0}$-module whose order is prime to $\textsf{char}(k_0)$. Recall that a cohomological invariant of $G$  over $k_0$ with coefficients in $C$ is a morphism of functors between $H^1(./k_0,G)$ and $H^*(./k_0,C)$. We denote by $\textsf{Inv}_{k_0}(G,C)$ the set of the cohomological invariants of $G$ over $k_0$ with coefficients in $C$.\\

Let $n\geq 2$. Let us recall that, for any field $k$, the pointed set $H^1(k,\mathfrak S_n)$ is in bijection with the set of isomorphism classes of \'etale $k$-algebras of rank $n$ (see for instance \cite{serre1984}).
In \cite{garibaldi2003}, Serre stated a splitting principle for cohomological invariants of the symmetric group $\mathfrak S_n$:
\begin{thm}[Serre, 2003]
\label{se1}
Let $k_0$ be a field such that $\textsf{char}(k_0)\neq 2$ and let $n\geq 2$. Let $a\in\textsf{Inv}_{k_0}(\mathfrak S_n,\mathbb Z/2\mathbb Z)$. Assume that, for every extension $k/k_0$, $a_k(E)=0$ whenever $E$ is an \'etale $k$-algebra isomorphic to a direct product of \'etale $k$-algebras of rank $\leq 2$. Then $a=0$.
\end{thm}

In this paper, we first generalize this result to finite Coxeter groups, provided that the base field has characteristic zero.\\

Let us recall the following definitions : 
\begin{defn}
Let $V$ be a $k_0$-vector space of finite dimension and let $r\in\textsf{End}_{k_0}(V)$. Then $r$ is a reflection of $V$ if the rank of $r-\textsf{id}$ is equal to $1$ and if $r^2=\textsf{id}$. 
\end{defn}
\begin{defn}
Let $W$ be a group of linear automorphisms of $V$. Then $W$ is a reflection group over $k_0$ if it is generated by reflections of $V$.
\end{defn}

Let us recall that finite reflection groups over $\mathbb R$ are exactly the finite Coxeter groups (see Section \ref{weylgroups} for further details).\\

Let $a\in \textsf{Inv}_{k_0}(W,C)$ and let $W'$ be a subgroup of $W$. Then the compositum of the morphisms of functors \begin{equation*}\xymatrix{H^1(./k_0,W')\ar[r] & H^1(./k_0,W) \ar[r]^{a} & H^*(./k_0,C)}\end{equation*} defines an invariant of $W'$, called the restriction of $a$ to $W'$ and denoted by $\textsf{Res}_W^{W'}(a)$. We then get a vanishing principle for finite Coxeter groups, already announced for Weyl groups by Serre (see 25.15, \cite{garibaldi2003}):
\begin{thm}
\label{iccox}
Let $W$ be a finite Coxeter group and let $k_0$ be a field of characteristic zero containing a subfield on which the real representation of $W$ as a reflection group is realizable. Let $C$ be a finite $\Gamma_{k_0}$-module and let $a\in \textsf{Inv}_{k_0}(W,C)$. Assume that every restriction of $a$ to an abelian subgroup of $W$ generated by reflections is zero. Then $a=0$.
\end{thm}


Note that the real representation of a Weyl group as a reflection group is realizable over $\mathbb Q$ (see Theorem \ref{weyl}), then Theorem \ref{iccox} is true for any Weyl group and any field $k_0$ of characteristic zero. In the case of the Weyl groups of type $A_n$, which means in the case of the symmetric groups, for any field $k$ of characteristic zero and for any maximal abelian subgroup $H$ generated by reflections (here transpositions), the image of $H^1(k,H)\rightarrow H^1(k,\mathfrak S_n)$ is exactly the set of the isomorphism classes of direct products of \'etale $k$-algebras of rank $\leq 2$. Therefore Theorem \ref{iccox} is exactly Theorem \ref{se1} in the case of the Weyl groups of type $A_n$.\\

As a first consequence of Theorem \ref{iccox}, we get that any normalized (i.e. vanishing on the trivial torsor) invariant $a\in \textsf{Inv}_{k_0}(W,C)$ is killed by $2$. We then give a characterization of the negligible cohomology classes of $W$.\\

Let $k_0$ be a field of characteristic different from $2$ and let $W$ be a Weyl group. Then the cup-product naturally endows the abelian group $\textsf{Inv}_{k_0}(W,\mathbb Z/2\mathbb Z)$ with the structure of an $H^*(k_0,\mathbb Z/2\mathbb Z)$-module. Using Theorem \ref{iccox}, we will show that it is free of finite rank for Weyl groups of classical types.\\

Let $k$ be a field of characteristic different from $2$ and let $m\geq 1$. Let us recall that any quadratic form $q$ of rank $m$ over $k$ may be written in an orthogonal basis $q\simeq\langle\alpha_1,...,\alpha_m\rangle$ for $\alpha_i\in k^\times$ (see for example \cite{serre1977}). For any $0\leq i\leq m$, if $q\simeq \langle\alpha_1,...,\alpha_m\rangle$, set \begin{equation*} w_i(q)=\underset{1\leq j_1<...<j_i\leq m}{\sum}(\alpha_{j_1})\cdot \ldots \cdot (\alpha_{j_i}).\end{equation*}  One may show (see \cite{delzant1962}) that, for $0\leq i\leq m$, $w_i(q)$ is well-defined and only depends on the isomorphism class of $q$. It yields cohomological invariants of the orthogonal group $\textbf{O}_n$ of the unit quadratic form of rank $n$, called Stiefel-Whitney invariants.\\

Recall that an \'etale $k$-algebra $L$ is a commutative $k$-algebra $L$ such that the $k$-bilinear form $(x,y)\mapsto \textsf{Tr}_{L/k}(xy)$ is non-degenerate, see for instance \cite{bourbaki2006}. Let $n\geq 2$. Recall also that $H^1(k,\mathfrak S_n)$ classifies the \'etale $k$-algebras of rank $n$ up to isomorphism. Let us now define some invariants of $\mathfrak S_n$. Set, for any $0\leq i\leq n$,  
\begin{equation*}\begin{aligned} w_i:H^1(./k_0,\mathfrak S_n) &\rightarrow H^*(./k_0,\mathbb Z/2\mathbb Z)\\ (L) &\mapsto w_i(q_L),\end{aligned}\end{equation*} where, for any \'etale $k$-algebra $L$, $q_L$ denotes the quadratic form $x\mapsto \textsf{Tr}_{L/k}(x^2)$. This invariant $w_i$ is called the $i^{th}$ Stiefel-Whitney invariant of $\mathfrak S_n$. We then have the following result (see \cite{garibaldi2003}, 25.13):
\begin{thm}[Serre, 2003] 
Let $k_0$ be a field of characteristic different from $2$ and let $n\geq 2$. Then the $H^*(k_0,\mathbb Z/2\mathbb  Z)$-module $\textsf{Inv}_{k_0}(\mathfrak S_n,\mathbb Z/2\mathbb Z)$ is free with basis $\{w_i\}_{0\leq i\leq [\frac{n}{2}]}$.
\end{thm}  

Let now $n\geq 2$ and let $W$ be a Weyl group of type $B_n$.
Let $k$ be a field of characteristic different from $2$. Then $H^1(k,W)$ classifies the pairs $(L,\alpha)$ up to isomorphism, where $L$ is an \'etale $k$-algebra and $\alpha$ a square-class in $L^\times$. We can now define two families of invariants of $W$.
Set, for any $0\leq i\leq n$,  
\begin{equation*}\begin{aligned} w_i:H^1(./k_0,\mathfrak S_n) &\rightarrow H^*(./k_0,\mathbb Z/2\mathbb Z)\\ (L,\alpha) &\mapsto w_i(q_L).\end{aligned}\end{equation*} 
Moreover, note that, for any $(L,\alpha)\in H^1(k, W)$, the isomorphism class of the quadratic form $q_{L,\alpha}:x\mapsto \textsf{Tr}_{L/k}(\alpha x^2)$ does not depend on the choice of a representative in the square-class $\alpha$ and is non-degenerate. 
Now set \begin{equation*}\begin{aligned} \widetilde{w}_i:H^1(./k_0,\mathfrak S_n) &\rightarrow H^*(./k_0,\mathbb Z/2\mathbb Z)\\ (L,\alpha) &\mapsto w_i(q_{L,\alpha}).\end{aligned}\end{equation*} 
These invariants are also called Stiefel-Whitney invariants of $W$. In this paper, we will prove the following result:
\begin{thm} 
\label{icbn}
Let $k_0$ be a field of characteristic zero, such that $-1$ and $2$ are squares in $k_0$. Let $n\geq 2$ and let $W$ be a Weyl group of type $B_n$. Then the $H^*(k_0,\mathbb Z/2\mathbb  Z)$-module $\textsf{Inv}_{k_0}(W,\mathbb Z/2\mathbb Z)$ is free with basis \begin{equation*}\{w_i\cdot\widetilde{w}_j\}_{\hspace{0.1cm} 0\leq i\leq [\frac{n}{2}],\hspace{0.2cm} 0\leq j\leq 2([\frac{n}{2}]-i)}.\end{equation*}
\end{thm}

Note that Weyl groups of type $C_n$ are isomorphic to Weyl groups of type $B_n$ (for any $n\geq 2$). \\

Let now $W$ be a Weyl group of type $D_n$ ($n\geq 4$). We have the exact sequence \begin{equation*} \xymatrix{1\ar[r] & W  \ar[r]& W' \ar[r]^{p} & \mathbb Z/2\mathbb Z\ar[r] & 1},\end{equation*} where $W'$ is a Weyl group of type $B_n$ and $p:(\epsilon_1,...,\epsilon_n,\sigma)\mapsto \overset{n}{\underset{i=1}{\prod}}\epsilon_i$.\\

Let $k_0$ be a field of characteristic different from $2$. If $a\in \textsf{Inv}_{k_0}(W',\mathbb Z/2\mathbb Z)$, then $\textsf{Res}_{W'}^W(a)$ is a cohomological invariant of $W$. Thus, for $0\leq i\leq n$, $\textsf{Res}_{W'}^W(w_i)$ is an invariant of $W$ and is still denoted by $w_i$. Likewise, for $0\leq i\leq n$, $\textsf{Res}_{W'}^W(\widetilde{w}_i)$ is an invariant of $W$ and is still denoted by $\widetilde{w}_i$. 
\begin{thm}
\label{icdn}
Let $k_0$ be a field of characteristic zero. Let $n\geq 4$ and let $W$ be a Weyl group of type $D_n$. Then the $H^*(k_0,\mathbb Z/2\mathbb  Z)$-module $\textsf{Inv}_{k_0}(W,\mathbb Z/2\mathbb Z)$ is free with basis \begin{equation*}\{w_i\cdot\widetilde{w}_j\}_{\hspace{0.1cm} 0\leq i\leq [\frac{n}{2}], \hspace{0.2cm} 0\leq j\leq 2([\frac{n}{2}]-i) \textrm{ and } j \textrm{ even }}.\end{equation*}
\end{thm}

\section*{Acknowledgments}
As explained in this introduction, the main goal of the following was to generalize Serre's results on cohomological invariants of the symmetric groups to Weyl groups, following the patterns of the proofs of Theorem 24.9 and Theorem 25.13 in \cite{garibaldi2003}, as suggested by Berhuy, our PhD advisor. As mentioned in \cite{garibaldi2003} 25.15, proofs of Theorem \ref{iccox} and the applications of Section \ref{app} were known to Serre. Notice also that Serre had already presented several applications of Theorem \ref{iccox} in talks and minicourses at Nottingham and Ascona conferences, especially to the alternating groups. Therefore, we warmly thank him to give his agreement with the publication of this paper. Eventually, we thank him for pointing out to us a mistake in an earlier version of this paper.\\

\section{Finite Coxeter groups and Weyl groups}
\label{weylgroups}
\begin{defn}
A Coxeter group $W$ is a group with a given presentation of type \begin{center} $\langle r_1,..,r_s\mid \forall i,j\in\{1,..,s\}, (r_ir_j)^{m_{i,j}}=1\rangle$,\end{center} where $\forall i,j\in\{1,...,s\}$, $m_{i,j}\in \mathbb N\cup \{+\infty\}$ and $m_{i,i}=1$ for every $i\in\{1,...,s\}$.
\end{defn}

It is well-known that a finite group $G$ is a Coxeter group if and only if it is a reflection group over $\mathbb R$ (see \cite{humphreys1990}). Note that there are some reflection groups over $\mathbb C$ which are not Coxeter groups (see for instance \cite{bourbaki1981},V.§5, exercise 4).
\begin{defn}
Let $V$ be a finite dimensional $\mathbb R$-vector space. A root system $S$ is a finite set of non-zero vectors in $V$ satisfying the conditions :
\begin{enumerate} 
\item for any $\alpha\in S$, $S\cap \mathbb R\alpha=\{\pm \alpha\}$
\item for any $\alpha\in S$, $r_{\alpha}(S)=S$, where $r_{\alpha}$ denotes the orthogonal reflection on  $V$ such that $\textsf{Im}(r_\alpha-\textsf{id}_V)=\mathbb R\alpha$.
\end{enumerate}
\end{defn}   

Note that a root system $S$ yields a finite Coxeter group (the group generated by the reflections $r_{\alpha}$, for any $\alpha\in S$). Conversely, any finite Coxeter group can be realized in this way, possibly for many different choices for $S$.\\

If a root system $S$ cannot be written $S_1\sqcup S_2$, with $S_1$ and $S_2$ two root systems, we say that $S$ is irreducible. Irreducible root systems are completely classified (and so are finite Coxeter groups) (see \cite{humphreys1990} or \cite{bourbaki1981} for details).\\

Let $(e_1,...,e_n)$ be a canonical basis of $\mathbb R^n$. Up to linear automorphism, irreducible root systems are classified in several types : 
\begin{itemize}
\item[] $A_n$ ($n\geq 1$) : let $V$ be the hyperplane of $\mathbb R^{n+1}$ such that the sum of coordinates equal to zero. Then $S=\{e_i-e_j\mid 1\leq i,j\leq n+1, i\neq j\}$. The Coxeter group is isomorphic to $\mathfrak S_{n+1}$.

\item[] $B_n$ (for $n\geq 2$) : $V=\mathbb R^n$, $S=\{\pm e_i, \pm e_j\pm e_l\mid 1\leq i\leq n, 1\leq j<l \leq n\}$. The Coxeter group is isomorphic to the semi-direct product $\big(\mathbb Z/2\mathbb Z\big)^n\rtimes \mathfrak S_n$, where $\mathfrak S_n$ acts on $(\mathbb Z/2\mathbb Z)^n$ by permuting coordinates.

\item[] $C_n$ (for $n\geq 2$) : $V=\mathbb R^n$, $S=\{\pm 2e_i, \pm e_j \pm e_l\mid 1\leq i\leq n, 1\leq j<l\leq n\}$. The Coxeter group is the same than in the type $B_n$.

\item[] $D_n$ ($n\geq 4$) : $V=\mathbb R^n$, $S=\{ \pm e_i \pm e_j)\mid 1\leq i<j\leq n\}$. The Coxeter group $W$ is defined by the exact sequence \begin{equation*} \xymatrix{1\ar[r] & W \ar[r]& W' \ar[r]^{p} & \mathbb Z/2\mathbb Z\ar[r] & 1},\end{equation*} where $W'$ is the reflection group of type $B_n$ and $\displaystyle p:(\epsilon_1,...,\epsilon_n,\sigma)\mapsto \prod_{i=1}^n\epsilon_i$. Moreover, $W$ is isomorphic to the semi-direct product $\big( \mathbb Z/2\mathbb Z\big)^{n-1}\rtimes \mathfrak S_n$.  

\item[] $E_6$ : $V=\{(x_i)_{1\leq i\leq 8}\in \mathbb R^8 \mid x_6=x_7=-x_8\}$, \begin{equation*}\begin{aligned}S=\{ \pm e_i \pm e_j, \pm\frac{1}{2}(e_8-e_7-e_6+\sum_{l=1}^5& (-1)^{\nu(l)}e_l) \\ &\mid 1\leq i<j\leq 5 \textrm{ and } \sum_{l=1}^5 \nu(l) \textrm{ even}\}\end{aligned}.\end{equation*} 

\item[] $E_7$ : let $V$ be the hyperplane of $\mathbb R^8$ orthogonal to $e_7+e_8$. Then  \begin{equation*} \begin{aligned} S=\{ \pm e_i \pm e_j, \pm(e_7-e_8) , \pm\frac{1}{2}(e_7 -e_8 &+\sum_{l=1}^6  (-1)^{\nu(l)}e_l) \\ &\mid 1\leq i<j\leq 6 \textrm{ and } \sum_{l=1}^6 \nu(l) \textrm{ odd}\}\end{aligned}.\end{equation*} 

\item[] $E_8$ : $V=\mathbb R^8$, \begin{equation*}S=\{ \pm e_i \pm e_j, \frac{1}{2}\sum_{l=1}^8 (-1)^{\nu(l)}e_l \mid 1\leq i<j\leq 8 \textrm{ and } \sum_{l=1}^8 \nu(l) \textrm{ even}\}.\end{equation*} 

\item[] $F_4$ :  $V=\mathbb R^4$, \begin{equation*}S=\{ \pm e_i, \pm e_j \pm e_l, \frac{1}{2}(\pm e_1 \pm e_2 \pm e_3 \pm e_4) \mid 1\leq i\leq 4, 1\leq j<l\leq 4 \}.\end{equation*} The Coxeter group is isomorphic to the semi-direct product \begin{equation*}\big(\big(\mathbb Z/2\mathbb Z\big)^3\rtimes \mathfrak S_4\big)\rtimes \mathfrak S_3,\end{equation*} where $ \big(\mathbb Z/2\mathbb Z\big)^3\rtimes \mathfrak S_4$ is the Coxeter group of type $D_4$ and $\mathfrak S_3$ acts on it by permuting vertices of the Dynkin diagram of $D_4$.

\item[] $G_2$ : let $V$ be the hyperplane of $\mathbb R^3$ with the sum of the coordinates equal to zero. Then \begin{equation*} \begin{aligned} S=\{ \pm (e_1-e_2), \pm (e_1-e_3), \pm (e_2-e_3), \pm (2e_1-e_2-e_3), &\pm (2e_2-e_1-e_3),\\ &\pm (2e_3-e_1-e_2)\}\end{aligned}.\end{equation*} The Coxeter group is isomorphic to the dihedral group $\mathbb D_6$ of order $12$.
\item[] $H_3$ : the Coxeter group is isomorphic to $\mathfrak A_5\times \mathbb Z/2\mathbb Z$, where $\mathfrak A_5$ denotes the alternating subgroup of $\mathfrak S_5$.
\item[] $H_4$ : the Coxeter group is the group of isometries of the hecatonicosahedroid.
\item[] $I_2(m)$, $m\geq 3$ : the Coxeter group is isomorphic to the dihedral group $\mathbb D_m$ of order $2m$.
\end{itemize}

With this classification, we get that every finite Coxeter group is isomorphic to a direct product of Coxeter groups of type $A$ to $I$.

\begin{defn}
A Weyl group $W$ is a finite Coxeter group with a root system $S$ satisfying the additional integrality condition : for any $\alpha,\beta\in S$, $ 2\frac{(\alpha,\beta)}{(\alpha,\alpha)}\in \mathbb Z$,\\ where $(.,.)$ denotes the usual scalar product.
\end{defn}

Any Weyl group is isomorphic to a direct product of groups of type $A$ to $G$. For these groups, we have the important following result (see \cite{springer1978}, Corollary 1.15) :

\begin{thm}
\label{weyl}
Let $W$ be a Weyl group. Every irreducible representation of $W$ is realizable over $\mathbb Q$. In particular, Weyl groups are reflection groups over $\mathbb Q$.
\end{thm}

Therefore, the real representation of a Weyl group as a real reflection group is realizable over $\mathbb Q$. By extension of scalars, Weyl groups are reflection groups over any field of characteristic zero. In particular, Theorem \ref{iccox} is true for any Weyl group and any field of characteristic zero.\\

Theorem \ref{weyl} is not true for a Coxeter group which is not a Weyl group. However, we have :
\begin{prop}
\label{cg} Let $W$ be a finite Coxeter group. There is a finite real extension $L$ of $\mathbb Q$ such that $W$ is a reflection group over $L$.
\end{prop}

Note that $L=\mathbb Q(\sqrt 5)$ for the Coxeter groups of type $H$ and $L=\mathbb Q(\cos(\frac{2\pi}{m}))$ for the Coxeter groups of type $I_2(m)$, for any $m\geq 3$ are the minimal fields such that Proposition \ref{cg} is satisfied.\\ 

Let $k_0$ be a field of characteristic zero. Thanks to the previous classification, $W$ is isomorphic to a direct product of groups of type $A$ to $I$. Then if $k_0$ contains the minimal field extensions over $\mathbb Q$ corresponding to the types in the decomposition of $W$, the representation of $W$ as a finite reflection group extends to $k_0$ and the assumption of Theorem \ref{iccox} is satisfied.\\

From now on, $W$ will denote a finite Coxeter group and $k_0$ a field  of characteristic zero containing a subfield on which the real representation of $W$ as a reflection group is realizable.

\section{ Vanishing theorem for cohomological invariants of finite Coxeter groups}
Let us recall the statement to prove :
\begin{thm*}
\label{iccox2}
Let $W$ be a finite Coxeter group and let $k_0$ be a field of characteristic zero containing a subfield on which the real representation of $W$ as a reflection group is realizable. Let $C$ be a finite $\Gamma_{k_0}$-module. Let also $a\in\textsf{Inv}_{k_0}(W,C)$. Assume that every restriction of $a$ to an abelian subgroup of $W$ generated by reflections is zero. Then $a=0$.
\end{thm*}

We will use the strategy suggested by Serre in \cite{garibaldi2003}, 25.15. Recall first that a cohomology class in $H^1(k, W)$ corresponds to an isomorphism class of a $W$-torsor over $k$. By \cite{garibaldi2003} Theorem 12.3, a cohomological invariant of $W$ is completely determined by its value on a versal torsor. Thanks to a Chevalley's theorem, we will construct a versal $W$-torsor $T^{\textsf{vers}}$ with rational base field $K=k_0(c_1,...,c_n)$. We will then show that, if a cohomological invariant $a$ of $W$ satisfies the hypothesis of Theorem \ref{iccox}, the cohomology class $a(T^{\textsf{vers}})$ is unramified at any place coming from an irreducible divisor of the affine space $\textsf{Spec}(k_0[c_1,...,c_n])$. By \cite{garibaldi2003}, Theorem 10.1, the cohomology class $a(T^{\textsf{vers}})$ is constant. Since $a$ vanishes on the trivial torsor, we get that $a=0$.

\subsection{Ramification of cohomology classes of $W$}
Let us recall what ramification means for cohomology classes in $H^1(k,W)$.
Let $R$ be a discrete valuation ring of valuation $v$, let $K$ be its fraction field and let $k$ be its residue field. 
Assume that $K$ is complete for the valuation $v$.
Let us also denote by $\Gamma_K$ (resp. by $\Gamma_k$) the absolute Galois group of $K$ (i.e. $\Gamma_K=\textsf{Gal}(K_{\textsf{sep}}/K)$, with $K_{\textsf{sep}}$ a fixed separable closure of $K$) (resp. the absolute Galois group of $k$). Finally, let us denote by $I_K$ the inertia group of $K$ and by $\pi:\Gamma_K\rightarrow \Gamma_k$ the quotient morphism.
\begin{prop}
\label{p1}
Let $\alpha\in H^1(K,W)$. If $\varphi$ is a cocycle representing $\alpha$, then the following assertions are equivalent :
\begin{itemize}
\item[(i)]$\varphi(I_K)=\{1_W\}$;
\item[(ii)] there is a unique group homomorphism $\overline{\varphi}:\Gamma_k\rightarrow W$ such that the following diagram is commutative : $\xymatrix{ \Gamma_K\ar[r]^{\varphi}\ar[d]_{\pi} & W\\ \Gamma_k\ar[ru]_{\overline{\varphi}} & }$;
\item[(iii)] $\alpha$ belongs to the image of the natural application $H^1(k,W)\rightarrow H^1(K,W)$.
\end{itemize} 
\end{prop}

Note that this statement only depends on the cohomology class $\alpha$.\\ 

\begin{dem}
\indent{   }
\begin{enumerate}
\item[\textit{(i)}$\Rightarrow$ \textit{(ii)}]Assume that $\varphi(I_K)=\{1_W\}$. Then $\varphi$ factors through  $\overline{\varphi}:\Gamma_K/I_K\rightarrow W$; yet $\Gamma_k=\Gamma_K/I_K$, so $\overline{\varphi}$ is the required morphism.
\item[\textit{(ii)}$\Rightarrow$ \textit{(iii)}] The homomorphism $\pi$ yields the map $\pi^*:H^1(k,W)\rightarrow H^1(K,W)$, given by $[\psi]\in H^1(k,W)\mapsto [\psi\circ\pi]\in H^1(K,W)$ (where $[.]$ denotes the cohomology class associated to the cocycle). Moreover, by \textit{(ii)}, since $\varphi=\overline{\varphi}\circ\pi$, $[\overline{\varphi}]$ is a preimage of $\alpha$ by $\pi^*$.  
\item[\textit{(iii)}$\Rightarrow$ \textit{(i)}] Assume that $\alpha$ admits a preimage $\beta\in H^1(k,W)$ by $\pi^*$. Then there is a cocycle $\psi$ representing $\beta$ such that $\varphi=\psi\circ \pi$, so the image of $I_K$ by $\varphi$ is trivial. $\blacksquare$ 
\end{enumerate}
\end{dem}

\begin{defn}
\label{ramcoh}
We say that the cohomology class $\alpha\in H^1(K,W)$ is unramified if $\alpha$ satisfies one of the three equivalent properties of Proposition \ref{p1}.\\ 
\end{defn}

\subsection{A versal $W$-torsor with rational base field}
\begin{defn}
\label{tv}
Let $G$ be a smooth algebraic group scheme over $k_0$. A versal $G$-torsor $P$ is a $G$-torsor over a finitely generated extension $K$ over $k_0$ such that there exists a smooth irreducible variety $X$ over $k_0$ with function field $K$ and a $G$-torsor $Q\rightarrow X$ with the following two properties : 
\begin{enumerate}
\item The fiber of $Q\rightarrow X$ at the generic point of $X$ is $P$;
\item For every extension $k$ over $k_0$, with $k$ infinite and for every $G$-torsor $T$ over $k$, the set $\{x\in X(k)\mid Q_x\simeq T\}$ is dense in $X$.
\end{enumerate} 
\end{defn}

As elements of $W$ are automorphisms of a vector space $V\simeq k_0^n$ for some $n>0$, $W$ naturally acts on the dual $V^*$ and on the associated symmetric algebra $\textsf{Sym}(V^*)$. Note that this $k_0$-algebra is isomorphic to a polynomial algebra  $k_0[x_1,...,x_n]$ with $n$ indeterminates. We then consider the underlying action of $W$ on $k_0[x_1,...,x_n]$ and the invariant subalgebra $k_0[x_1,...,x_n]^W$. By a theorem of Chevalley (\cite{bourbaki1981}, 5.5.5, p. 115), $k_0[x_1,...,x_n]^W$ is a polynomial $k_0$-algebra of transcendence degree $n$. In other words, $k_0[x_1,...,x_n]^W \simeq k_0[c_1,...,c_n]$ for some independent indeterminates $c_1,...,c_n$ over $k_0$.\\

We now translate it into scheme language. Set $Q=\textsf{Spec}(k_0[x_1,...,x_n])$ and $ X=\textsf{Spec}(k_0[c_1,...,c_n])$. We have a morphism $f:Q\rightarrow X$ which is exactly the quotient morphism $\textsf{Aff}^n_x\rightarrow \textsf{Aff}^n_x/W=\textsf{Aff}^n_c$. Let $y$ be an element of $k_0[x_1,...,x_n]$ whose orbit by $W$ has maximal order. We then localize $f$ at the locus $\Delta_c=f(\Delta_x)$, where $\Delta_x=\{w.y-w'.y\mid w\neq w', w,w'\in W\}$. We then get from $f$ a morphism $Q_{\Delta_x}\rightarrow X_{\Delta_c}$ that we still denote by $f$. With this localization, $W$ acts without fixed points on $Q_{\Delta_x}$ and we still have $Q_{\Delta_x}/W=X_{\Delta_c}$. Hence, $Q_{\Delta_x}$ is a $W$-torsor with base $X_{\Delta_c}$.\\

We denote by $K=k_0(c_1,...,c_n)$ the function field of $X$ (which is also the function field of $X_{\Delta_c}$) and by $L=k_0(x_1,...,x_n)$ the function field of $Q$. As $X_{\Delta_c}$ is an irreducible variety, let us denote by $T^{\textsf{vers}}$ the fiber of $f$ at the unique generic point of $X_{\Delta_c}$. Thus, $T^{\textsf{vers}}$ is a $W$-torsor over $K$, corresponding to the field extension $L/K$ which is Galois, with Galois group $W$.

\begin{prop}
\label{pvers}
Keeping the notation above, $T^{\textsf{vers}}$ is a versal torsor  for $W$ over $k_0$.
\end{prop}

\begin{dem}
Let $k/k_0$ be a field extension. Let $T$ be a $W$-torsor over $k$. Then $T$ corresponds to a Galois $W$-algebra over $k$ and we choose a generator $(a_1,...,a_n)$. We localize $Q_{\Delta_x}$ at the ideal $(x_1-a_1,...,x_n-a_n)$ of $k[x_1,...,x_n]_{\Delta_x}$.  The image $x=f(x_1-a_1,...,x_n-a_n)$ is a $k$-point of $X_{\Delta_c}$ and the fiber of $f$ in $x$ is a $W$-torsor over $k$ isomorphic to $T$. As $k$ is infinite (since $k_0$ has characteristic zero), the set of the generators of the Galois $W$-algebra is dense with respect to the Zariski topology on $\textsf{Aff}_x^n$, so condition \textit{2.} in Definition \ref{tv} is satisfied. $\blacksquare$
\end{dem}

Note that the isomorphism class of $T^{\textsf{vers}}$ corresponds to the cohomology class of the natural projection $\varphi^{\textsf{vers}}:\begin{aligned}\Gamma_K &\rightarrow W \\ \gamma &\mapsto \gamma_{\mid L}\end{aligned}$.

\subsection{Ramification of the versal torsor $T^{\textsf{vers}}$}
In this section, we want to study the ramification of the isomorphism class of the versal torsor $T^{\textsf{vers}}$ at the different valuations on $K$ which are trivial on $k_0$.
These valuations are determined by the irreducible divisors of $\textsf{Aff}^n_c$. 
Let $D$ be such a divisor. Let us denote by $v_D$ the discrete valuation on $K$ associated to $D$, $K_D$ the completion of $K$ with respect to this valuation and $k_0(D)$ the residue field of $K$ for $v_D$, which identifies with the function field of $D$ over $k_0$.
We denote by $T^{\textsf{vers}}_D$ the image of $T^{\textsf{vers}}$ under the application \begin{equation*} \begin{aligned} H^1(K,W)&\rightarrow H^1(K_D,W)\\ [\varphi] &\mapsto [\varphi\circ i_D]\end{aligned}\end{equation*} where $i_D:\Gamma_{K_D}\rightarrow \Gamma_K$ is the natural inclusion.\\

The aim of this paragraph is to study the ramification of the cohomology class of $T^{\textsf{vers}}_D$. We denote by $\varphi^{\textsf{vers}}_D$ the morphism $\xymatrix{\Gamma_{K_D}\ar[r]^{i_D} & \Gamma_K \ar[r]^{\varphi^{\textsf{vers}}} &W}$; it represents the cohomology class of $T^{\textsf{vers}}_D$.
Thus, by Proposition \ref{p1}, we have to study the subgroup $\varphi^{\textsf{vers}}_D(I_{K_D})$.\\

Since $[\varphi^{\textsf{vers}}]$ is represented by the Galois extension $L/K$, $[\varphi^{\textsf{vers}}_D]$ is represented by the Galois $W$-algebra $L\otimes_KK_D$ over $K_D$. By \cite{neukirch1999} II.8, there is an isomorphism of $K_D$-algebras \begin{equation*} L\otimes_KK_D\simeq \underset{\widetilde{v}\mid v_D} {\prod} L_{\widetilde v}\end{equation*} where $L_{\widetilde v}$ denotes the completion of $L$ with respect to the extension $\widetilde{v}$ of $v_D$.\\

Let $\widetilde{v_D}$ be an extension of the valuation $v_D$ to $L$. We denote by $\widetilde{ L_D}$ the completion of $L$ with respect to this valuation. Then $\widetilde{L_D}$ is a Galois extension of $K_D$, with Galois group $\widetilde W =\{w\in W, \widetilde{v_D}\circ w=\widetilde{v_D}\}$, which is of course a subgroup of $W$.\\ 

Set $e=(0,...,0,1,0,...,0)$ in the product $\underset{\widetilde{v}\mid v_D} {\prod} L_{\widetilde v}$, where $1\in \widetilde{L_D}$. Then $e$ is a primitive idempotent of $L\otimes_KK_D$ and by \cite{knus1998}, Proposition 18.18, $\widetilde{L_D}=e.(L\otimes_KK_D)$ is a Galois $\widetilde W$-algebra (and a field) and we have the isomorphism of $W$-algebras \begin{equation*}L\otimes_KK_D\simeq \textsf{Ind}_{\widetilde W}^W(\widetilde{L_D}).\end{equation*}

Thus, since the induced algebra (for the Galois algebras) corresponds to the inclusion for the cocycles, $\varphi^{\textsf{vers}}_D$ factors through $\widetilde W$ : \begin{center} $\xymatrix{\Gamma_{K_D}\ar[r]^{\varphi^{\textsf{vers}}_D} \ar[d]_{\psi}& W\\ \widetilde W\ar@{^{(}->}[ru] & }$,\end{center} where $\psi$ is a cocycle representing the cohomology class corresponding to $\widetilde{L_D}/K_D$. It yields that $\varphi^{\textsf{vers}}_D(I_{K_D})=\psi(I_{K_D})$. Therefore the ramification of $T^{\textsf{vers}}_D$ is $\psi(I_{K_D})$.\\

Let us denote by $l(\widetilde{v_D})$ the residue field associated with $\widetilde{L_D}$ and if $w\in \widetilde W$, let us denote by $\overline w$ the induced $k_0(D)$-automorphism (as $w$ respects the valuation $\widetilde{v_D}$, $w$ restricts to $w:\mathcal O_{\widetilde{v_D}}\rightarrow \mathcal O_{\widetilde{v_D}}$, where $\mathcal O_{\widetilde{v_D}}$ denotes the valuation ring of $\widetilde{v_D}$ in $L$ and sends the maximal ideal of $\mathcal O_{\widetilde{v_D}}$ in itself, so going to quotients, we get an automorphism $\overline w$ of $l(\widetilde{v_D})$). We then introduce the inertia subgroup $\widetilde{I}=\{w\in\widetilde W\mid \overline w=\textsf{id}_{l(\widetilde{v_D})}\}$ of $\widetilde W$. 
\begin{lem}
\label{psi}
Keeping the notation above, $\psi(I_{K_D})\subset\widetilde I$.
\end{lem}

\begin{dem}
Let us denote by $\overline{k_0(D)}$ an algebraic closure of $k_0(D)$. Recall that $(K_D)_{\textsf{sep}}$ has residue field $\overline{k_0(D)}$ and that  $(k_0(D))_{\textsf{sep}}$ is the residue field corresponding to the biggest subextension of $(K_D)_{\textsf{sep}}$ fixed by the inertia group $I_{K_D}$. Let $\gamma\in I_{K_D}$. Then the $k_0(D)$-automorphism $\overline{\gamma}$ is trivial over $k_0(D)_{\textsf{sep}}$. In other words, the image of $\gamma$ by the group homomorphism $\Gamma_{K_D}\rightarrow\Gamma_{k_0(D)}$ is the identity and we have the commutative diagram \begin{center}$\xymatrix{\Gamma_{K_D}\ar[r] \ar[d]_{\psi}& \Gamma_{k_0(D)}\ar[d]^{\textsf{res}}\\ \widetilde W\ar[r] & \textsf{Gal}(l(\widetilde{v_D})/k_0(D))}$\end{center} where horizontal maps are induced by going to quotients (by valuation theory, the sequence
\begin{equation*}\xymatrix{0\ar[r] & \widetilde I\ar[r] & \widetilde W\ar[r] & \textsf{Gal}(l(\widetilde{v_D})/k_0(D)))\ar[r] & 0}\end{equation*} is exact). Then the $k_0(D)$-automorphism of $l(\widetilde{v_D})$ induced by $\psi(\gamma)$ is equal to the identity, which  proves that $\psi(\gamma)$ belongs to $\widetilde I$. $\blacksquare$
\end{dem}

Let us now study the inertia group $\widetilde I$ and let us introduce the discriminant $\textsf{Discr}(L/K)$ of $L/K$. Let us recall that the isomorphism class of the versal torsor $T^{\textsf{vers}}$ identifies with the isomorphism class of the Galois algebra $L/K$, i.e. with the set of $K$-embeddings of $L$ in $K_{\textsf{sep}}$. Yet, these embeddings are completely determined by the image of a primitive element $y$ of $L$ over $K$.
Therefore, this discriminant may be written as:
\begin{center}
$\textsf{Discr}(L/K)=\underset{ t\neq t'}{\prod}(t(y)-t'(y))$,
\end{center}
where $t,t':L\hookrightarrow K_{\textsf{sep}}$.\\

Moreover, one can choose, as a primitive element, a polynomial in $k_0[x_1,...,x_n]$ with total degree 1 : as $k_0$ is infinite, there exists $y= a_1x_1+...+a_nx_n$ (where $a_i\in k_0$ for $i=1,...,n$) such that, for all $w\neq w'\in W$, $w(y)\neq w'(y)$. Indeed, as $W$ is a group, it is enough to check that there exist $a_1,...,a_n\in k_0$ such that $y=a_1x_1+...+a_nx_n$ and that, for all $w\in W$, $w(y)\neq y$, which is satisfied as soon as the vector $\begin{pmatrix}a_1\\ \vdots \\a_n\end{pmatrix}$ is not an eigenvector of any matrix representing a non-trivial element of $W$ in the basis $(x_1,...,x_n)$ of the dual space $V^*\simeq k_0^n$.\\

From now on, $y$ will denote a  primitive element of $L/K$, which is a polynomial of total degree 1 in $x_1,...,x_n$. Let us now compute the ramification of $T^{\textsf{vers}}_D$.
\begin{lem}
\label{p3}
Assume that $D$ is an irreducible divisor which does not divide the ideal $\big(\textsf{Discr}_{L/K}\big)$.
Then the isomorphism class of $T^{\textsf{vers}}_D$ is unramified.
\end{lem} 

\begin{dem}
Since we have shown above that the ramification is contained in $\widetilde I$, it is enough to prove that $\widetilde I$ is trivial. Yet the sequence
\begin{center}$\xymatrix{0\ar[r] & \widetilde I\ar[r] & \widetilde W\ar[r] & \textsf{Gal}(l(\widetilde{v_D})/k_0(D)))\ar[r] & 0}$\end{center}
is exact. As $D$ does not divide the discriminant, the extension $\widetilde{L_D}/K_D$ is unramified, which shows that $[\widetilde{L_D}:K_D]=[l(\widetilde{v_D}):k_0(D)]$, so $\widetilde W\simeq \textsf{Gal}(l(\widetilde{v_D})/k_0(D)))$. Therefore, $\widetilde I$ is trivial. $\blacksquare$
\end{dem}

\begin{lem}
\label{p4}
Assume now that $D$ is an irreducible divisor of $\textsf{Spec}(k_0[c_1,...,c_n])$ which divides the discriminant ideal $\big(\textsf{Discr}_{L/K}\big)$. Then $\widetilde I=\langle r \rangle$, where $r$ is a reflection of $W$.
\end{lem}

\begin{dem}
Since $D$ divides the discriminant ideal, the extension $\widetilde{L_D}/K_D$ is ramified, so its inertia group $\widetilde I$ is not trivial. Let us compute it.
Since $\widetilde{v_D}$ is a valuation on $L=k_0(x_1,...,x_n)$, which is trivial on $k_0$ (because it extends $v_D$ which is itself trivial on $k_0$), $\widetilde{v_D}$ is a valuation associated with an irreducible divisor $\widetilde D$ of $\textsf{Spec}(k_0[x_1,...,x_n])$; furthermore, since $\widetilde{v_D}$ extends $v_D$, the divisor $\widetilde D$ is above $D$, so $\widetilde D$ is generated by an irreducible factor of $\textsf{Discr}(L/K)$ decomposed in $k_0[x_1,...,x_n]$. Then there are two distinct elements $t_1,t_2\in T^{\textsf{vers}}$ such that $\widetilde D=(t_1(y)-t_2(y))$  (we identify the image of $L$ in $K_{\textsf{sep}}$ with $L$ itself; that is why we consider $t_1(y)$ and $t_2(y)$ as polynomials in $x_1,...,x_n$). Therefore, the valuation $\widetilde{v_D}$ is described as follows : for any $f\in L$, $\widetilde{v_D}(f)$ is equal to the order of $(t_1(y)-t_2(y))$ as zero or pole in the rational fraction $f$.\\

Let $w\in \widetilde I$. Then $\overline w=\textsf{id}_{l(\widetilde{v_D})}$. Let $f\in \mathcal O_{\widetilde{v_D}}$ (i.e. which has not $t_1(y)-t_2(y)$ as a pole). As $l(\widetilde{v_D})=\mathcal O_{\widetilde{v_D}}/\mathcal M_{\widetilde {v_D}}$, there is a $g\in\mathcal O_{\widetilde{v_D}}$ such that :
\begin{equation*}
w(f)=f+g.(t_1(y)-t_2(y))\end{equation*} 
If we now write $g=\frac{g_0}{g_1}$, with $g_0,g_1\in k_0[x_1,...,x_n]$ and $t_1(y)-t_2(y)$ not dividing $g_1$, we get that \begin{equation*} g_1.(w(f)-f)=g_0.(t_1(y)-t_2(y)).\end{equation*}

Consider the particular case where $f$ is a polynomial in $k_0[x_1,...,x_n]$. Then the equality now reads in 
$k_0[x_1,...,x_n]$ (because $W$ acts on $k_0[x_1,...,x_n]$) and since $t_1(y)-t_2(y)$ does not divide $g_1$, $t_1(y)-t_2(y)$ divides $w(f)-f$. Therefore, there is a polynomial $g_2$ in $k_0[x_1,...,x_n]$ such that \begin{equation}
\label{wf}
w(f)-f=g_2.(t_1(y)-t_2(y))
\end{equation}

Assume now that $f$ is a polynomial of total degree 1. Then $f$ identifies with a linear form on $V$ and as $w$ is an automorphism of $V$, $w(f)=f\circ w^{-1}$ still is a linear form on $V$, so via the identification $\textsf{Sym}(V^*)\simeq k_0[x_1,...,x_n]$, $w(f)$ still is a polynomial of total degree 1. \\

Let us now show that the total degree of $t_1(y)-t_2(y)$ is equal to $1$.
First note that, as $t_1(y)-t_2(y)$ is an irreducible factor of $\textsf{Discr}_{L/K}$ in $k_0[x_1,...,x_n]$, it has total degree $\geq 1$. Assume that $w\neq \textsf{id}$. Then there is a $f_0\in V^*$ such that $w(f_0)\neq f_0$. Thus, $w(f_0)-f_0$ is a polynomial of total degree $0$ or $1$. Taking total degree in Equation (\ref{wf}), we get that  $t_1(y)-t_2(y)$ has exactly total degree $1$. Therefore, $g_2$ has total degree at most $0$. Thus, for any $f\in V^*$, there exists $a\in k_0$ such that  \begin{center}$w(f)-f=a.(t_1(y)-t_2(y))$.\end{center}

We then get that $w$ is a pseudo-reflection of $V^*$ (i.e. an endomorphism such that the rank of $w-\textsf{id}_{V^*}$ is equal to 1). Yet, since $W$ is a reflection group over $\mathbb R$, the only pseudo-reflections in $W$ are reflections. Therefore, the non-trivial elements of $\widetilde I$ are reflections.\\
 
By \cite{serre1979} IV.2, Corollary 2 of Proposition 7, we get that $\widetilde I$ is cyclic (note that the residue field of $l(\widetilde{v_D})$ is an extension of $k_0$, then has characteristic zero). Then $\widetilde I$ is of order 2 (recall that it can not be trivial). Finally $\widetilde I=\{1,r\}$ (where $r$ is a reflection of $W$). $\blacksquare$
\end{dem}

Let us recall that we have $\varphi^{\textsf{vers}}_D(I_{K_D})=\psi(I_{K_D})\subset \widetilde I=\{1,r\}$. Then \begin{center} $\varphi^{\textsf{vers}}_D(I_{K_D})=\{1\}$ or $\varphi^{\textsf{vers}}_D(I_{K_D})=\{1,r\}$.\end{center} In the first case, $T^{\textsf{vers}}$ is unramified in $D$. In the second case, we state the key lemma for our inductive proof of Theorem \ref{iccox} :
\begin{lem}
\label{kl}
Assume that $\varphi^{\textsf{vers}}_D(I_{K_D})=\langle r \rangle$. 
Then there is a subgroup $W_0$ of $W$ generated by reflections, such that $\varphi^{\textsf{vers}}_D(\Gamma_{K_D})\subset W_0\times \langle r \rangle \subset W$.
\end{lem}

\begin{dem}
Since the sequence \begin{equation*}\xymatrix{1 \ar[r] & I_{K_D} \ar[r] & \Gamma_{K_D} \ar[r] & \Gamma_{k_0(D)}\ar[r] & 1}\end{equation*} is exact, $\varphi^{\textsf{vers}}_D(I_{K_D})$ is normal in $\varphi^{\textsf{vers}}_D(\Gamma_{K_D})$. Therefore, $r$ is in the center of $\varphi^{\textsf{vers}}_D(\Gamma_{K_D})$, that is to say that $\varphi^{\textsf{vers}}_D(\Gamma_{K_D})$ is contained in the centralizer $C(r)$ of $r$ in $W$.\\
 
By assumption on $k_0$, the real representation $W\hookrightarrow GL(V_{\mathbb R})$ of $W$ as a reflection group over $\mathbb R$ yields a representation $W\hookrightarrow GL(V_{k_0})$ of $W$ as a reflection group over $k_0$. 

Let now $e$ be a non-zero vector of $\textsf{Im}(r-\textsf{id}_{V_{\mathbb R}})$ and let $H$ be the hyperplane of the fixed points of $r$ in $V_{\mathbb R}$. Let also $w\in W$. Then $w$ and $r$ commute if and only if $\mathbb Re$ and $H$ are stable by $w$ (see for instance \cite{bourbaki1981}, Proposition 3 of p. 68). 
Assume that $w$ and $r$ commute. Then, $w(H)\subset H$ and $w(e)=b.e$ for some $b\in \mathbb R$. As $W$ is finite, $b$ is a root of the unity in $\mathbb R$, so $b=\pm 1$. \\

Let $W_0=\{w\in W \mid w(e)=e\}$. As an isotropy subgroup of $W$, $W_0$ is a reflection group over $\mathbb R$ (see \cite{humphreys1990}, Theorem of p. 22) and hence still remain a reflection group over $k_0$.\\ 

It remains to prove that $C(r)\simeq W_0\times \langle r \rangle$. Let us first note that, as $W$ acts by isometries on the euclidean space $V_\mathbb R$, $W_0=\{w\in W\mid w(H)\subset H \textrm{ and } w(e)=e\}$. One can show easily that $W_0$ and $\langle r \rangle$ are normal in $C(r)$, that the intersection $W_0\cap \langle r\rangle$ is trivial and that $W_0.\langle r \rangle =C(r)$. Therefore, $C(r)$ is the direct product of  $W$ and $\langle r \rangle$.  $\blacksquare$
\end{dem}

\subsection{Proof of Theorem \ref{iccox}}

Let us begin with a key lemma for our proof by induction:
\begin{lem}
\label{p6}
Let $W$ be a finite Coxeter group and let $k_0$ be a field satisfying the hypothesis of Theorem \ref{iccox}. Let $W'$ be a proper subgroup of $W$ which is also generated by reflections. Let us assume that there is a reflection $r$ of $W$ which is not in $W'$ and which commutes with any reflection of $W'$. If Theorem \ref{iccox} is true for $W'$, then Theorem \ref{iccox} is also true for $W'\times \langle r\rangle$.
\end{lem}

\begin{dem}
Let $a\in\textsf{Inv}_ {k_0}(W'\times \langle r\rangle, C)$ such that every restriction to an abelian subgroup generated by reflections is zero. Let $k/k_0$ be a field extension. We have the isomorphism $H^1(k,W'\times \langle r \rangle)\simeq H^1(k,W')\times H^1(k,\langle r \rangle)$. Then, in the sequel of the proof, we will denote the elements of $H^1(k,W' \times \langle r \rangle)$ by pairs $(\alpha,\epsilon)$, where $\alpha$ is a cohomology class in $H^1(k,W')$ and $\epsilon$ a square-class in $H^1(k,\langle r \rangle)$.\\

Let $(\alpha_0,\epsilon_0)\in H^1(k,W'\times \langle r \rangle)$ be such an element. For any extension $k'/k$, we set \begin{equation*}\begin{aligned}(\widetilde a_{\epsilon_0,k})_{k'}:H^1(k',W')&\rightarrow H^*(k',C)\\ \alpha &\mapsto a_{k'}(\alpha,\epsilon_0)\end{aligned}.\end{equation*} It is easily seen that these maps define a cohomological invariant $\widetilde a_{\epsilon_0,k}$ of $W'$ over $k$. Assume that $\alpha\in H^1(k',W')$ lies in the image of a map $H^1(k',H')\rightarrow H^1(k', W')$, where $H'$ is an abelian subgroup of $W'$ generated by reflections. Then $(\alpha,\epsilon_0)$ is in the image of $H^1(k',H'\times \langle r \rangle)\rightarrow H^1(k', W'\times \langle r \rangle)$ and since $H'\times \langle r\rangle $ is an abelian subgroup of $W'\times \langle r \rangle$ generated by reflections, by assumption on $a$, $a_{k'}(\alpha,\epsilon_0)=0$. Hence, $(\widetilde a_{\epsilon_0,k})_{k'}(\alpha)=0$. Therefore, $\widetilde a_{\epsilon_0, k}$ satisfies the assumption of Theorem \ref{iccox}. As Theorem \ref{iccox} is true for $W'$ (by the hypothesis on $W'$), we get that $\widetilde a_{\epsilon_0,k}=0$. It then yields that $a_k(\alpha_0,\epsilon_0)=0$. Finally, $a=0$. $\blacksquare$
\end{dem}

We can now give the proof of Theorem \ref{iccox}.
\begin{dem}
For convenience, we say that $a\in \textsf{Inv}_{k_0}(W,C)$ satisfies $(P)$ if every restriction of $a$ to an abelian subgroup generated by reflections is zero.
We will show Theorem \ref{iccox} by induction on the order $m$ of $W$ : if $m=1$ or $m=2$, it is trivial. Let $m\geq 3$. Assume that, for every integer $l$ with $1\leq l<m$, every cohomological invariant of a Coxeter group which satisfies the assumption of Theorem \ref{iccox} over $k_0$ of order $l$ satisfying $(P)$ is zero.\\ 

Let $a\in\textsf{Inv}_{k_0}(W,C)$ satisfying $(P)$. We will prove that, for any irreducible divisor $D$ of $\textsf{Aff}^n_c$, the residue $r_{v_D}(a_K(\varphi^{\textsf{vers}}))$, at the valuation $v_D$ corresponding to the divisor $D$, is zero. Then, by Theorem 10.1 in \cite{garibaldi2003}, $a_K(\varphi^{\textsf{vers}})$ will be constant and since  $\varphi^{\textsf{vers}}$ corresponds to the versal torsor  $T^{\textsf{vers}}$ for $W$ over $k_0$, $a$ will be constant by \cite{garibaldi2003}, Theorem 12.3. Thus, since the restrictions of $a$ to any abelian subgroup generated by reflections are zero,  $a$ will vanish on the trivial torsor and we will get that $a=0$.\\ 

Let $D$ be an irreducible divisor in $\textsf{Spec}(k_0[c_1,...,c_n])$. Let us prove that the residue 
$r_{v_D}(a_K(\varphi^{\textsf{vers}}))$ is zero. We know by Lemma \ref{p3} and Lemma \ref{p4} that 
$\varphi^{\textsf{vers}}_D(I_{K_D})=\{1\}$ or $\varphi^{\textsf{vers}}_D(I_{K_D})=\langle r \rangle$ for some reflection $r\in W$. In the first case, $\varphi^{\textsf{vers}}_D$ is not ramified so, by \cite{garibaldi2003} Theorem 11.7, $r_{v_D}(a_K(\varphi^{\textsf{vers}}))=r_{v_D}(a_{K_D}(\varphi^{\textsf{vers}}_D))=0$.\\

Assume now that $\varphi^{\textsf{vers}}_D(I_{K_D})=\langle r\rangle$. By Lemma \ref{kl}, $\varphi^{\textsf{vers}}_D(\Gamma_{K_D})\subset W_0\times \langle r \rangle$. Since $W_0$ is a proper subgroup of $W$ and a reflection group over $k_0$, it satisfies the assumptions of Theorem \ref{iccox} so by the induction hypothesis, Theorem \ref{iccox} is true for $W_0$. By Lemma \ref{p6}, Theorem \ref{iccox} is also true for the group $W_0\times \langle r\rangle$. \\

Since $a$ satisfies $(P)$, $\textsf{Res}_{\hspace{0.4cm} W }^{W_0\times \langle r\rangle}(a)$ also satisfies $(P)$, so  $\textsf{Res}_{\hspace{0.4cm} W }^{W_0\times \langle r\rangle}(a)=0$. Thus, as $[\varphi^{\textsf{vers}}_D]$ lies in the image of the map $H^1(K_D,W_0\times \langle r\rangle)\rightarrow H^1(K_D, W)$, we get that $a_{K_D}(\varphi^{\textsf{vers}}_D)=0$. Hence, its residue $r_{v_D}(a_{K_D}(\varphi^{\textsf{vers}}_D))$ is also zero.\\

We then have shown that, for every irreducible divisor of $\textsf{Spec}(k_0[c_1,...,c_n])$,\\ $r_{v_D}(a_{K_D}(\varphi^{\textsf{vers}}_D))=0$. This concludes the proof. $\blacksquare$
\end{dem}


\section{Applications}
\label{app}
The following two applications directly generalize similar results of Serre for the symmetric groups and were already known to Serre (see \cite{garibaldi2003} 25.15).
\subsection{ Invariants are killed by $2$}
Recall that $W$ is a finite Coxeter group, $k_0$ is a field of characteristic zero containing a subfield on which the representation of $W$, as a reflection group is realizable. Recall also that $C$ is a finite $\Gamma_{k_0}$-module.
\begin{defn}
An element $a\in\textsf{Inv}_{k_0}(W,C)$ is normalized if $a_{k_0}([1])=0$, where $[1]$ denotes the cohomology class of the trivial cocycle in $H^1(k_0,W)$.
\end{defn}

Let us state a first consequence of Theorem \ref{iccox} (generalizing \cite{garibaldi2003}, 24.12).
\begin{cor}
\label{2k}
For every normalized invariant $a\in\textsf{Inv}_{k_0}(W,C)$, $2a=0$. In particular, if $C$ has odd order, there is no non-trivial normalized invariant.
\end{cor}

\begin{dem}
Let $a\in \textsf{Inv}_{k_0}(W,C)$ be  normalized. By Theorem \ref{iccox}, it is enough to prove that, for any abelian subgroup $H$ of $W$ generated by reflections, $2\textsf{Res}_W^H(a)=0$. We prove it by induction on the order $m$ of $W$. For $m=1$ or $2$, it is trivial.\\ 

Let $m\geq 3$ and let $W$ be a reflection group over $k_0$ of order $m$. Let denote by $S$ a root system corresponding to $W$. Let $H$ be an abelian subgroup of $W$ generated by reflections. Then $H\simeq \langle r_1\rangle\times \cdots \times \langle r_s\rangle$ for some $s\geq 1$ and for some pairwise commuting reflections $r_1,...,r_s$ in $W$. Let $e_1$ be the root in $S$  corresponding to $r_1$ and let $W'$ be the group generated by reflections given by the root subsystem $S'=S\cap \{e_1\}^{\perp}$.  
Then $W'$ is a proper subgroup of $W$ and a reflection group over $k_0$. Using the induction hypothesis with $W'$ and with the normalized invariant $\textsf{Res}_W^{W'}(a)$, we get that \begin{center}$2(\textsf{Res}_W^{W'}(a))=0$.\end{center}

Let $H'=\langle r_2\rangle\times \cdots \times \langle r_s\rangle$. Since $H'\subset W'$, \begin{equation*}2\textsf{Res}_W^{H'}(a)=2\textsf{Res}_{W'}^{H'}(\textsf{Res}_W^{W'}(a))=0.\end{equation*}
Let $k$ be an extension of $k_0$ and let $T\in H^1(k,H)$. Then $T=T_1\times T_2$, where $T_1\in H^1(k,\langle r_1\rangle)$ and $T_2\in H^1(k,H')$. Thus, \begin{equation*}2\textsf{Res}_W^{H'}(a)_k(T_2)=0.\end{equation*}  
Now set $T'=T_1'\times T_2$ the $H$-torsor , where $T_1'$ is the trivial torsor in $H^1(k,\langle r\rangle)$. Since $T'=T_1'\times T_2=i^*(T_2)$ with $i:H'\hookrightarrow H$ and $i^*:H^1(k,H')\rightarrow H^1(k,H)$ the induced map, the definition of the restriction map yields \begin{equation*}(\textsf{Res}_{W}^{H}(a))_k(T')= \textsf{Res}_H^{H'}(\textsf{Res}_{W}^{H}(a))_k(T_2).\end{equation*} Therefore, \begin{center} $2\textsf{Res}_W^{H}(a)_k(T')=0$.\end{center}

Moreover, there is an extension $k'/k$ of degree at most $2$ such that $T_1$ and $T_1'$ are isomorphic over $k'$. Then $T'$ and $T$ also are isomorphic over $k'$. Hence, \begin{equation*}\textsf{Res}_W^H(a)_{k'}(\textsf{Res}_{k'/k}(T))=\textsf{Res}_W^H(a)_{k'}(\textsf{Res}_{k'/k}(T')),\end{equation*} so \begin{equation*}\textsf{Res}_{k'/k}(\textsf{Res}_W^H(a)_{k}(T))=\textsf{Res}_{k'/k}(\textsf{Res}_W^H(a)_{k}(T'))\end{equation*} and applying the corestriction map $\textsf{Cor}_{k'/k}$, we get that  \begin{equation*}[k':k].\textsf{Res}_W^H(a)_k(T)=[k':k].\textsf{Res}_W^H(a)_k(T')\end{equation*} which proves that
 \begin{equation*}2\textsf{Res}_W^H(a)_k(T)=2\textsf{Res}_W^H(a)_k(T')=0.  \end{equation*} We conclude by using Theorem \ref{iccox} to $2a$. $\blacksquare$
\end{dem}

Corollary \ref{2k} allows us to restrict to $\Gamma_{k_0}$-modules $C$ of even order. The most elementary one is of course $\mathbb Z/2\mathbb Z$ endowed with the trivial action of $\Gamma_{k_0}$ and we will take $C=\mathbb Z/2\mathbb Z$ in most of our examples.

\subsection{Application to negligible cohomology}
Let $G$ be a finite group and let $M$ a finite $G$-module, with trivial action. Let $k$ be a field and  $x\in H^*(G,M)$. We have a natural map \begin{equation*}\begin{aligned}(a_x)_k: H^1(k,G)&\rightarrow H^*(k,M)\\ [\varphi] &\mapsto [\varphi^*(x)]\end{aligned}\end{equation*}
Since $G$ acts trivially on $M$, for any extension $k'/k$, the maps $(a_x)_{k'}$ define a cohomological invariant of $G$ over $k$ with coefficients in $M$. \\

This gives us a family of cohomological invariants of $G$. We want to determine the cohomology classes $x\in H^*(G,M)$ for which these invariants are zero. 

\begin{defn}
Let $G$ be a finite group and let $M$ be a $G$-module. Then a cohomology class $x\in H^*(G,M)$ is negligible if, for any field $k$ and any (continuous) homomorphism $\varphi:\Gamma_k\rightarrow G$, we have $\varphi^*(x)=0$ in $H^*(k,M)$.
We denote by $H^*_{\textsf{negl}}(G,M)$ the subset of $H^*(G,M)$ consisting of the negligible cohomology classes.
\end{defn}

In fact, it is enough to consider only fields of characteristic zero (see \cite{garibaldi2003}, 26.1):
\begin{prop}
\label{n1}
An element $x\in H^*(G,M)$ is negligible if $\varphi^*(x)=0$ for any field $k$ of characteristic zero and any $\varphi:\Gamma_k\rightarrow G$.
\end{prop}

For any $x\in H^*(G,M)$, let us denote by $a_x$ the cohomological invariant over $\mathbb Q$ induced by $x$. Then Proposition \ref{n1} exactly says that \begin{equation*} H^*_{\textsf{negl}}(G,M)=\{x\in H^*(G,M)\mid a_x=0\}.\end{equation*}

As a first example, let us give a negligibility criterion for $2$-elementary groups (see \cite{garibaldi2003}, Lemma 26.4).
\begin{ex}
\label{n2}
Let $G$ be a $2$-elementary group and let $x\in H^i(G,\mathbb Z/2\mathbb Z)$. Then $x$ is negligible if and only if the restriction of $x$ to every subgroup  of $G$ of order $\leq 2$ is zero. 
\end{ex}

The following result is clear and its proof is left to the reader.
\begin{lem}
\label{n3}
Let $i\geq 0$ and let $x\in H^i(G,M)$. For any subgroup $H$ of $G$, \begin{equation*}\textsf{Res}_G^H(a_x)=a_{\textsf{Res}_G^H(x)}.\end{equation*}
\end{lem}

The next result is the natural generalization of a result of Serre on negligible cohomology classes of $\mathfrak S_n$ (see \cite{garibaldi2003}, 26.3) to Weyl groups.
\begin{thm}
Let $W$ be a Weyl group and let $M$ be a finite $W$-module, with trivial action. Let $i\geq 0$.
We have the following assertions :
\begin{enumerate}
\item[(1)] $x\in H^i(W,M)$ is negligible if and only if its restrictions to the abelian subgroups generated by reflections  are negligible.
\item[(2)] for any $i>0$, for any $x\in H^i(W,M)$, the cohomology class $2x$ is negligible.
\item[(3)] An element $x\in H^i(W,\mathbb Z/2\mathbb Z)$ is negligible if, and only if, its restrictions to the subgroups of order $\leq 2$ of $W$ are zero.
\end{enumerate}
\end{thm}

\begin{dem}
\begin{enumerate}
\item[\textit{(1)}] Let $x\in H^i(W,M)$.  By Proposition \ref{n1}, the kernel of the natural map \begin{equation*}\begin{aligned} H^i(W,M) &\rightarrow \textsf{Inv}_{\mathbb Q}(W,M)\\ x & \mapsto a_x\end{aligned}\end{equation*} is exactly $H^i_{\textsf{negl}}(W,M)$. 
Furthermore, since $W$ is a reflection group over $\mathbb Q$ (see Theorem \ref{weyl}), Theorem \ref{iccox} yields that $a_x=0$ if and only if $\textsf{Res}_W^H(a_x)=0$ for any abelian subgroup $H$  of $W$ generated by reflections. Thus, by Lemma \ref{n3}, $x$ is negligible if and only if for any abelian subgroup $H$ of $W$ generated by reflections, the restriction $\textsf{Res}_W^H(x)$ is negligible.
  
\item[\textit{(2)}] Let $i>0$ and $x\in H^i(W,M)$.  Let us show that $a_{2x}=2a_x$. For any field $k$ of characteristic zero and any $[\varphi]\in H^1(k,W)$, $(a_{2x})_k([\varphi])$ is represented by  $(\gamma_1,...,\gamma_i)\mapsto 2x (\varphi(\gamma_1),...,\varphi(\gamma_i)) $. Hence, $a_{2x}=2a_x$. We conclude the proof by using Corollary \ref{2k}.

\item[\textit{(3)}] By \textit{(1)}, $x\in H^i(W,\mathbb Z/2\mathbb Z)$ is negligible if and only if for any abelian subgroup $H$ of $W$ generated by reflections, $\textsf{Res}_W^H(x)$ is negligible. Then, Example \ref{n2} allows us to conclude. $\blacksquare$
\end{enumerate}
\end{dem}

\section{Cohomological invariants of Weyl groups of type $B_n$}
Let $k_0$ be a field of characteristic different from $2$ and let $W$ be a Weyl group. Then the cup-product endows the abelian group $\textsf{Inv}_{k_0}(W,\mathbb Z/2\mathbb Z)$ with an $H^*(k_0,\mathbb Z/2\mathbb Z)$-module structure. 

\subsection{Restriction to a $2$-elementary subgroup}
Let us first give the full description of the cohomological invariants of a $2$-elementary group with coefficients in $\mathbb Z/2\mathbb Z$ (see \cite{garibaldi2003}, 16.4).
\begin{ex}
\label{ex1}
Let $k_0$ be a field of characteristic different from $2$ and let $n\geq 1$. Let $H=(\mathbb Z/2\mathbb Z)^n$.  Recall that $H^1(k, H)\simeq H^1(k,\mathbb Z/2\mathbb Z)\times \cdots \times H^1(k,\mathbb Z/2\mathbb Z)$ and that $H^1(k,\mathbb Z/2\mathbb Z)\simeq k^\times/k^{\times 2}$. Let  $I\subset \{1,...,n\}$. For any $k/k_0$, let us define  \begin{equation*}\begin{aligned}(a_I)_k:H^1(k,H) &\rightarrow H^{i}(k,\mathbb Z/2\mathbb Z)\\ (x_1,...,x_n) &\mapsto (x)_I\end{aligned}\end{equation*} where $i$ denotes the cardinality of $I$ and $(x)_I$ denotes the cup product of the $(x_i)$, $i\in I$. It is clear that these maps define a cohomological invariant of $H$. We then have (see \cite{garibaldi2003}, Theorem 16.4):
\begin{prop}
\label{e1}
The $H^*(k_0,\mathbb Z/2\mathbb Z)$-module $\textsf{Inv}_{k_0}(H,\mathbb Z/2\mathbb Z)$ is free with basis $\{a_I\}$, where $I$ describes the subsets of $\{1,...,n\}$.
\end{prop}
\end{ex}

Let us come back to the general case. Let $H$ be a subgroup of $G$. Since the restriction of $a\in \textsf{Inv}_{k_0}(G,C)$ is the compositum \begin{equation*}\xymatrix{H^1(./k_0,H)\ar[r] & H^1(./k_0,G)\ar[r]^a & H^*(./k_0,C)},\end{equation*} it defines a map \begin{equation*}\xymatrix{\textsf{Inv}_{k_0}(G,C)\ar[r]^{\textsf{Res}_G^{H}} & \textsf{Inv}_{k_0}(H,C)}.\end{equation*} Let $N_{G}(H)$ be the normalizer of $H$ in $G$. For any $a\in \textsf{Inv}_{k_0}(H,C)$, for any $g\in N_{G}(H)$ and for any cocycle $\varphi:\Gamma_k\rightarrow H$ on $k/k_0$, \begin{equation*}(g.a)_k([\varphi])=[g\varphi g^{-1}]\end{equation*} (recall that $[.]$ denotes the associated cohomology class). Thus, we get an action of $N_G(H)$ on $\textsf{Inv}_{k_0}(H,C)$ and by \cite{garibaldi2003}, Proposition 13.2, the image of the restriction map $\textsf{Res}_G^{H}$ is contained in $\textsf{Inv}_{k_0}(H,C)^{N_G(H)/H}$ (it is clear that the action of $H$ is trivial).

\begin{ex}
\label{ex2}
Let $W$ be a Weyl group of type $B_n$ ($n\geq 2$). Its associated root system is $S=\{\pm e_i,\pm e_i\pm e_j, 1\leq i\neq j\leq n\}$. Set $S_0=\{\pm e_{1},...,\pm e_n\}$. The subgroup $H_0$ of $W$ generated by the reflections corresponding to the roots in $S_0$ is clearly isomorphic to $(\mathbb Z/2\mathbb Z)^n$.
By Example \ref{ex1}, the family of invariants $\{a_I\}_{I\subset \{1,...,n\}}$ is a basis of the $H^*(k_0,\mathbb Z/2\mathbb Z)$-module $\textsf{Inv}_{k_0}(H_0,\mathbb Z/2\mathbb Z)$. Let us identify the invariants fixed by the normalizer $N_0$ of $H_0$. Since $W$ permutes the lines $\mathbb Re_i$ ($1\leq i\leq n$), we have $N_0=W$. Moreover, we have the exact split sequence \begin{equation*}\xymatrix{1 \ar[r] & H_0 \ar[r] & W \ar[r] & \mathfrak S_n\ar[r] & 1}.\end{equation*} Thus, $N_0/H_0\simeq \mathfrak S_n$ and acts on $H_0$ by permuting the coordinates. 
\begin{prop}
\label{e2}
For $0\leq i\leq n$, the cohomological invariants $a_i^{(0)}=\underset{I\subset \{1,...,n\}; \mid I\mid=i}{\sum} a_I$ form a basis of the submodule $\textsf{Inv}_{k_0}(H_0,\mathbb Z/2\mathbb Z)^{N_0/H_0}$.
\end{prop}
\end{ex}

By using Proposition \ref{e1}, the proof of Proposition \ref{e2} is easy and left to the reader.

\subsection{The vanishing principle for Weyl groups of type $B_n$}
\label{s1}
Let $n\geq 2$, let $(e_1,...,e_n)$ be the canonical basis of $\mathbb R^n$ and let $S$ be the root system of type $B_n$ : $S=\{\pm e_i,\pm e_i\pm e_j, 1\leq i\neq j\leq n\}$. Let us denote by $W$ its Weyl group. By the classification given in Section \ref{weylgroups}, $W$ is isomorphic to the semi-direct product $\big(\mathbb Z/2\mathbb Z\big)^n\rtimes \mathfrak S_n$, where $\mathfrak S_n$ acts on $\big(\mathbb Z/2\mathbb Z\big)^n$ by permuting coordinates.\\

Let us first give an interpretation of the first Galois cohomology set of $W$. Let $k$ be a field of characteristic different from $2$. We call pointed \'etale $k$-algebra of rank $n$ any couple $(L,\alpha)$ with $L$ an \'etale $k$-algebra of rank $n$ and $\alpha$ a square-class in $L^\times$. Let $L,L'$ be \'etale $k$-algebras of rank $n$ and let $\alpha, \alpha'$ be square-classes respectively in $L^\times$ and in $L'^\times$. A morphism of pointed \'etale algebras $(L,\alpha)\rightarrow (L',\alpha')$ is a homomorphism $f:L[\sqrt{\alpha}] \rightarrow L'[\sqrt{\alpha'}]$ of $k$-algebras such that $f(L)\subset L'$ and $f(\sqrt{\alpha})=\lambda \sqrt{ \alpha'}$ for some $\lambda \in L'$ (note that by $L[\sqrt{\alpha}]$, we mean the $k$-algebra $k[X]/(X^2-\alpha)$).

\begin{prop}
\label{h1bn}
Let $k$ be any field of characteristic different from $2$. The set $H^1(k, W)$ is in bijection with the set of isomorphism classes of the pointed \'etale $k$-algebras of rank $n$. Moreover the trivial cocycle  is mapped onto the base point $(k^n,1)$. 
\end{prop}

\begin{dem}
By \cite{knus1998}, 29.12, as the pairs $(L,\alpha)$ are exactly the twisted $k$-forms of $(k^n,1)$, up to isomorphism, we just have to see that the automorphism group of $(k^n,1)$ is $W$. Note that $k^n[\sqrt 1]=k^{2n}$ and that $\textsf{Aut}_k(k^{2n})=\mathfrak S_{2n}$. Let us call $1,2,...,2n$ the $2n$ factors of $k^{2n}$ and let us consider the elements of $\textsf{Aut}_k(k^n,1)$ as permutations of $\{1,...,2n\}$. We have the morphism \begin{equation*}\begin{aligned}\Phi :\textsf{Aut}_k(k^n,1)&\rightarrow \textsf{Aut}_k{k^n}\\ f & \mapsto f_{\mid k^n}. \end{aligned}\end{equation*} Then it is easily seen that $\textsf{Ker}(\Phi)$ consists of the permutations fixing the subsets $\{1,2\},\{3,4\},...,\{2n-1,2n\}$. Thus $\textsf{Ker}(\Phi)$ is isomorphic to $(\mathbb Z/2\mathbb Z)^n$. It yields the exact sequence \begin{equation*}\xymatrix{1\ar[r] & (\mathbb Z/2\mathbb Z)^n \ar[r] & \textsf{Aut}_k(k^n,1) \ar[r]^-{\Phi} & \mathfrak S_n \ar[r] & 1  }.\end{equation*} Moreover, the inclusion $\mathfrak S_n=\textsf{Aut}_k(k^n)\hookrightarrow W$  splits this exact sequence and $\mathfrak S_n$ acts on $\textsf{Ker}(\Phi)$ by permuting the subsets $\{1,2\},\{3,4\},...,\{2n-1,2n\}$. Therefore $\textsf{Aut}_k(k^n,1)\simeq W$. $\blacksquare$
\end{dem}

More explicitely, the isomorphism class of $(L,\alpha)$ is represented, as a cohomology class, by the cocycle \begin{equation*}
\varphi_{L,\alpha}:\begin{aligned}\Gamma_k &\rightarrow W\\ \gamma &\mapsto \big( (\epsilon_1(\gamma), ..., \epsilon_n(\gamma)), \sigma_\gamma\big) \end{aligned}
\end{equation*} where, for any $\gamma\in \Gamma_k$:
\begin{itemize}
\item[] $\sigma_\gamma$ is the permutation induced by the action of $\gamma$ on $X(L)=\textsf{Hom}(L,k_{\textsf{sep}})$;
\item[] $\epsilon_i(\gamma)=1$ if $\gamma$ does not exchange factors $2i-1$ and $2i$ in $L[\sqrt \alpha]\otimes_kk_{\textsf{sep}}\simeq k_{\textsf{sep}}^{2n}$, $\epsilon_i(\gamma)=-1$ otherwise.\\
\end{itemize}

For any integer $q$ such that $0\leq q\leq [\frac{n}{2}]$, let $H_q$ be the subgroup of $W$ associated with the root subsystem of $S$ :\begin{equation*}S_q=\{\pm e_1\pm e_2,\pm e_3\pm e_4,...,\pm e_{2q-1}\pm e_{2q},\pm e_{2q+1},...,\pm e_n\}.\end{equation*} Then it is easily seen that the set $\{H_q\mid 0\leq q\leq [\frac{n}{2}]\}$ form a system of representatives modulo conjugation of the maximal abelian subgroups of $W$ generated by reflections. Assume that $n$ is even (the case $n$ odd is similar), and let $k/k_0$ be an extension. Let $0\leq q\leq \frac{n}{2}$ and let $u_1,...,u_{n/2},v_1,...,v_{n/2}$ be square-classes in $k^\times$. The image of $(u_1,v_1,u_2,v_2,...,u_{n/2},v_{n/2})$ by the map $H^1(k,H_q)\rightarrow H^1(k, W)$ is \begin{equation*}T_q=(k(\sqrt{u_1v_1})\times \cdots \times k(\sqrt{u_qv_q})\times k^{n-2q}, (u_1,u_2,...,u_q,u_{q+1},v_{q+1},...,u_{n/2},v_{n/2})). \end{equation*}

We may now reformulate Theorem \ref{iccox} for Weyl groups of type $B_n$, $n\geq 2$.
\begin{cor}
\label{icbn1}
Let $k_0$ be a field of characteristic zero, let $C$ be a finite $\Gamma_{k_0}$-module and let $a\in \textsf{Inv}_{k_0}(W,C)$. Then $a=0$ if, and only if, for any $0\leq q\leq [\frac{n}{2}]$, $\textsf{Res}_W^{H_q}(a)=0$. In other words, $a=0$ if and only if $a$ vanishes on the pairs \begin{equation*}(k(\sqrt{t_1})\times \cdots\times k(\sqrt{t_q})\times k^{n-2q},(\alpha_1,...,\alpha_{n-q}))\end{equation*} (for $0\leq q\leq [\frac{n}{2}]$), where the square-class $\alpha_i$ has a representative in $k^\times$ for any $0\leq i \leq n-q$.
\end{cor}

\begin{dem}
Since $\Gamma_k$ acts trivially on $W$, the images of the maps $H^1(k,H)\rightarrow H^1(k,W)$ and   $H^1(k,H')\rightarrow H^1(k, W)$ are the same if $H$ and $H'$ are two conjugate subgroups of $W$. Then the result directly follows from Theorem \ref{iccox}. 
$\blacksquare$
\end{dem}

\subsection{Cohomological invariants of $\mathbb D_4$}
The results of this subsection were presented by Serre in his minicourse at the Ascona conference in 2007 in a different way. Recall that the Weyl group of type $B_2$ is isomorphic to the dihedral group $\mathbb D_4$, with $8$ elements. We will directly determine (i.e. without using Theorem \ref{iccox}) the cohomological invariants of $\mathbb D_4$, by computing residues in $H^*(./k_0,\mathbb Z/2\mathbb Z)$. For details about Galois cohomology with coefficients in $\mathbb Z/2\mathbb Z$, see for instance \cite{berhuy2010}. \\

If $k=k_0(c_1,...,c_r)$ is a rational field extension over $k_0$ with transcendance degree $r$ and if $P$ is an irreducible polynomial (in the variables $c_1,...,c_r$ over $k_0$), let us denote by $D_P$ the irreducible divisor in $\textsf{Spec}(k_0[c_1,...,c_r])$ associated with $P$. Let us also denote by $v_P$ the valuation $v_{D_P}$ corresponding to the divisor $D_P$ and by $r_P$ the residue map $r_{v_P}$.\\

Let us now state a technical lemma :
\begin{lem}
\label{prel1}
Let $l\geq 0$ and let $k=k_0(t,u,v_1,...,v_l)$ be a rational field extension over $k_0$ with transcendence degree $l+2$. If $\alpha\in H^*(k,\mathbb Z/2\mathbb Z)$ is not ramified at any $k_0$-valuation on $k$, except maybe at the valuations $v_t$ and $v_u$, then there exist $c_0,c_1,c_2,c_3\in H^*(k_0,\mathbb Z/2\mathbb Z)$ such that \begin{equation*}\alpha= c_0+c_1\cdot (t) +c_2\cdot (u) +c_3\cdot (u) \cdot (t).\end{equation*}
\end{lem}

\begin{dem}
First set $k'=k_0(t,u)$. If $l\geq 1$, then, for any effective divisor $D$ of $\textsf{Spec}(k'[v_1,...,v_l])$, $\alpha$ is not ramified at $v_D$. By \cite{garibaldi2003}, Theorem 10.1, $\alpha\in H^*(k',\mathbb Z/2\mathbb Z)$. If now $D$ is an irreducible divisor of $\textsf{Spec}(k_0(t)[u])$, different from $D_u$, then $\alpha$ is not ramified at $v_D$. By \cite{garibaldi2003} 9.4, there exist some $\alpha_0,\alpha_1\in H^*(k_0(t),\mathbb Z/2\mathbb Z)$ such that \begin{equation}\label{eq25} \alpha=\alpha_0+\alpha_1\cdot (u).\end{equation} 

Let now $D$ be an irreducible divisor of $\textsf{Spec}(k_0[t])$ different from $D_t$. Then $D\times \mathbb A^1$ is an irreducible divisor of $\textsf{Spec}(k_0[u,t])$ different from $D_t$ and $D_u$ and we have the commutative diagram (see \cite{garibaldi2003} 8.2)
\begin{equation*}\xymatrix{ H^*(k_0(t),\mathbb Z/2\mathbb Z) \ar[r]^{r_{v_D}} \ar[d]_{\textsf{Res}_{k_0(t,u)/k_0(t)}} & H^*(\kappa(D),\mathbb Z/2\mathbb Z) \ar[d]^{\textsf{Res}_{\kappa(D)(u)/\kappa(D)}} \\ H^*(k_0(t,u),\mathbb Z/2\mathbb Z) \ar[r]_{r_{v_{D\times \mathbb A^1}}}&  H^*(\kappa(D)(u),\mathbb Z/2\mathbb Z),}\end{equation*} 
which implies that \begin{equation*} 0=r_{v_{D\times \mathbb A^1}}(\alpha)= \textsf{Res}_{\kappa(D)(u)/\kappa(D)}(r_{v_D}(\alpha_0))+\textsf{Res}_{\kappa(D)(u)/\kappa(D)}(r_{v_D}(\alpha_1))\cdot (u).\end{equation*}  Since $\kappa(D)(u)/\kappa(D)$ is purely transcendental, the map $\textsf{Res}_{\kappa(D)(u)/\kappa(D)}$ is injective, so \begin{equation*} 0=r_D(\alpha_0)+r_D(\alpha_1)\cdot (u),\end{equation*} and since $u$ is an indeterminate over $\kappa(D)$, \begin{equation*}r_D(\alpha_0)=r_D(\alpha_1)=0.\end{equation*} 

Therefore, $\alpha_0$ and $\alpha_1$ are not  ramified at any $k_0$-valuation on $k_0(t)$ except maybe at $v_t$, so, by \cite{garibaldi2003} 9.4, there exist some $c_0,c_1,c_2,c_3 \in H^*(k_0,\mathbb Z/2\mathbb Z)$ such that \begin{center}$\alpha_0=c_0+c_1\cdot (t)$ and $\alpha_1=c_2+c_3\cdot (t)$,\end{center} and Equation (\ref{eq25}) allows us to conclude.
$\blacksquare$\\
\end{dem}

\begin{thm}
\label{d4}
Let $k_0$ be a field of characteristic different from $2$. Then the set $\textsf{Inv}_{k_0}(\mathbb D_4,\mathbb Z/2\mathbb Z)$ is a free $H^*(k_0,\mathbb Z/2\mathbb Z)$-module with basis $\{1,w_1,\widetilde{w}_1,\widetilde{w}_2\}$.
\end{thm}

\begin{dem}
Let $K=k_0(t,u,v)$, where $t,u,v$ are independent indeterminates. We consider the $\mathbb D_4$-torsor $\big( K(\sqrt t),u+v\sqrt t \big)$ over $K$. This torsor is versal over $k_0$ for $\mathbb D_4$ (see Definition \ref{tv}). Recall that, by \cite{garibaldi2003}, Theorem 12.3, every cohomological invariant is completely determined by its value on a versal torsor. Moreover, we have the following formulae :
\begin{equation*}\begin{aligned} w_1(K(\sqrt t),u+v\sqrt t)&=(2.2t)=(t),\\ 
\widetilde {w}_1(K(\sqrt t),u+v\sqrt t)&=(2u. 2ut(u^2-v^2t))=(t(u^2-v^2t)) \textrm{ and}\\ 
\widetilde {w}_2(K(\sqrt t),u+v\sqrt t)&=(2u)\cdot (2ut(u^2-v^2t))=(2u)\cdot (-t(u^2-v^2t)).\end{aligned}\end{equation*}

It then remains to prove the following two facts :
\begin{enumerate}
\item[(i)] the family $\{1,(t), (t(u^2-v^2t)), (2u)\cdot (-t(u^2-v^2t))\}$ is free in the module $H^*(K,\mathbb Z/2\mathbb Z)$ over $H^*(k_0,\mathbb Z/2\mathbb Z)$;
\item[(ii)] if $a\in \textsf{Inv}_{k_0}(\mathbb D_4,\mathbb Z/2\mathbb Z)$, there exist $d_0,d_1,\widetilde{d}_1,\widetilde{d}_2\in H^*(k_0,\mathbb Z/2\mathbb Z)$ such that \begin{equation*}a_K(K(\sqrt t),u+v\sqrt t) = d_0 + d_1\cdot (t) + \widetilde{d}_1\cdot (t(u^2-v^2t)) + \widetilde{d}_2\cdot (2u)\cdot (-t(u^2-v^2t)).\end{equation*}
\end{enumerate}

\noindent Let us first show (i).
Let $\lambda_0,\lambda_1,\lambda_2,\lambda_3\in H^*(k_0,\mathbb Z/2\mathbb Z)$ such that \begin{equation*} 0=\lambda_0+\lambda_1\cdot (t)+\lambda_2\cdot (t(u^2-v^2t))+\lambda_3\cdot (2u)\cdot (-t(u^2-v^2t)).\end{equation*} 
Let us take the residue at the valuation corresponding to $(u^2-v^2t)$. Then, \begin{equation*} 0=\lambda_2+\lambda_3\cdot (2u)=(\lambda_2+\lambda_3\cdot (2))+ \lambda_3\cdot (u).\end{equation*}
Taking now the residue at $v_u$, it is easily seen that $\lambda_3=0$, which implies that $\lambda_2=0$. We then obtain that $0=\lambda_0+\lambda_1\cdot (t)$, so taking the residue at $v_t$, we get that $\lambda_1=0$. Thus, $\lambda_0=0$ and this proves (i).\\ 

Let us prove (ii). Let $a$ be a cohomological invariant of $\mathbb D_4$ over $k_0$ and set \begin{equation*} \beta=a_K(K(\sqrt t),u+v\sqrt t)\end{equation*} and \begin{equation*}\beta_1=\textsf{Res}_{K(\sqrt t)/K}(\beta)=a_{K(\sqrt t)}\big( K(\sqrt t)^2, (u+v\sqrt t,u-v\sqrt t) \big).\end{equation*} Let us consider $H_0\subset \mathbb D_4$ (viewed as the Weyl group of type $B_2$) introduced in Example \ref{ex2}. Since the cohomology class associated with the pointed algebra $\big( K(\sqrt t)^2, (u+v\sqrt t,u-v\sqrt t) \big)$ lies in the image of $H^1(K(\sqrt{t}),H_0)\rightarrow H^1(K(\sqrt t),\mathbb D_4)$, we get that \begin{equation*}\beta_1={\textsf{Res}_W^{H_0}(a)}_{K(\sqrt t)} \big( u+v\sqrt t,u-v\sqrt t \big).\end{equation*} 

Therefore, since $\textsf{Res}^{H_0}_W(a)$ is a cohomological invariant of $H_0\simeq \big( \mathbb Z/2\mathbb Z \big)^2$, by Example \ref{ex2}, there are some $b_0,b_1,b_2\in H^*(k_0,\mathbb Z/2\mathbb Z)$ such that, for any field extension $k/k_0$ and for any $(\alpha_1,\alpha_2)\in H^1(k,H_0)$, we have : \begin{equation*}(\textsf{Res}_W^{H_0}(a))_k(\alpha_1,\alpha_2)=b_0+b_1\cdot (\alpha_1\alpha_2)+b_2\cdot (\alpha_1)\cdot(\alpha_2).\end{equation*} 

Hence, we get that \begin{equation*}\beta_1=b_0+b_1\cdot (u^2-v^2t)+b_2\cdot (u+v\sqrt t)\cdot (u-v\sqrt t).\end{equation*}
Since the extension $K(\sqrt t)/K$ is not ramified at the valuation corresponding to $(u^2-v^2t)$, we have the commutative diagram (\cite{garibaldi2003}, Proposition 8.2) :
\begin{equation*}\xymatrix{ H^*(K,\mathbb Z/2\mathbb Z) \ar[r]^{r_{u^2-v^2t}} \ar[d]_{\textsf{Res}_{K(\sqrt t)/K}} & H^*(k_0(u,v),\mathbb Z/2\mathbb Z) \ar[d]^{\textsf{id}} \\ H^*(K(\sqrt t),\mathbb Z/2\mathbb Z) \ar[r]_{r_{u+v\sqrt t}}&  H^*(k_0(u,v),\mathbb Z/2\mathbb Z)}.\end{equation*} 
In the residue field associated with $r_{u+v\sqrt t}$ over $K(\sqrt t)$, we have $\overline{2u}=\overline{u-v\sqrt t}$, since $u+v\sqrt t+u-v\sqrt t=2u$. The previous commutative diagram yields that \begin{equation*}r_{u^2-v^2t}(\beta)=r_{u+v\sqrt t}(\beta_1)=b_1+b_2\cdot (2u).\end{equation*}  In particular, $r_{u^2-v^2t}(\beta)$ is not ramified, except maybe at $v_u$. Now set 
\begin{equation*}\beta'=\beta+r_{u^2-v^2t}(\beta)\cdot (u^2-v^2t).\end{equation*} Let us show that the cohomology class $\beta'$ is not ramified except maybe at $v_t$ and at $v_u$. Let $D$ be an irreducible divisor of $\textsf{Spec}(k_0[t,u,v])$ different from $D_t$, $D_u$ and $D_{u^2-v^2t}$. With Definition \ref{ramcoh}, it is easily seen that the cohomology class of the pointed algebra $(K(\sqrt t), u+v\sqrt t)$ is not ramified except maybe at the valuations corresponding to the irreducible divisors of the discriminant of the algebra $K(\sqrt t)\big(\sqrt{u+v\sqrt t}\big)/K$, i.e. at $v_t$ and at $v_{u^2-v^2t}$. Since $D\neq D_t,D_{u^2-v^2t}$, the cohomology class of the pointed algebra $(K(\sqrt t), u+v\sqrt t)$ is not ramified at $v_D$, thus $r_{v_D}(\beta)=0$. Hence, \begin{equation*}r_{v_D}(\beta')=r_{v_D}(\beta)+r_{v_D}((b_1+b_2\cdot (2u))\cdot (u^2-v^2t))=0.\end{equation*} Furthermore, \begin{equation*} r_{u^2-v^2t}(\beta')=r_{u^2-v^2t}(\beta)+r_{u^2-v^2t}\big(r_{u^2-v^2t}(\beta)\cdot (u^2-v^2t)\big)=2r_{u^2-v^2t}(\beta)=0.\end{equation*} 

Hence $\beta'$ is not ramified, except maybe at $v_t$ and $v_u$. By Lemma \ref{prel1}, there exist $c_0,c_1,c_2,c_3\in H^*(k_0,\mathbb Z/2\mathbb Z)$, such that \begin{equation*} \beta'=c_0+c_1\cdot (t)+c_2\cdot (u)+c_3\cdot (u)\cdot (t), \end{equation*} so 
\begin{equation*}\beta =c_0+c_1\cdot (t)+c_2\cdot (u)+c_3\cdot (u)\cdot (t)+(b_1+b_2\cdot (2u))\cdot (u^2-v^2t).\end{equation*}

Yet we know that $r_u(\beta)=0$, therefore : 
\begin{equation*}
0=c_2+c_3\cdot (t)+b_2\cdot (\overline{u^2-v^2t})=c_2+c_3\cdot (t)+b_2\cdot (-t).
\end{equation*}

It yields that \begin{equation*}0=c_2+b_2\cdot (-1) +(c_3+b_2)\cdot (t).\end{equation*} As the family $\{1,(t)\}$ is free in the $H^*(k_0,\mathbb Z/2\mathbb Z)$-module $H^*(k_0(t),\mathbb Z/2\mathbb Z)$, we get that $c_2=b_2\cdot (-1)$ and $c_3=b_2$. Eventually, 
\begin{equation*} \begin{aligned}
\beta & = c_0 +c_1\cdot (t)+b_1\cdot (u^2-v^2t) + b_2\cdot (2)\cdot (u^2-v^2t)+ b_2\cdot (u)\cdot (-t(u^2-v^2t)) \\& = c_0 +(c_1+b_2\cdot (2))\cdot (t)+b_1\cdot (u^2-v^2t) + b_2\cdot (2u)\cdot (-t(u^2-v^2t)).
\end{aligned}\end{equation*}
Hence,  
\begin{equation*}\beta = c_0+(c_1+b_1+b_2\cdot (2))\cdot (t)+b_1\cdot (t(u^2-v^2t)) + b_2\cdot (2u)\cdot (-t(u^2-v^2t)),\end{equation*}  This proves (ii).
$\blacksquare$
\end{dem}

Let us end this section by giving relations between some cohomological invariants of $\mathbb D_4$ that we will use below :
\begin{prop}
\label{reld4}
We have the following equalities : \begin{center}$w_2=(2)\cdot w_1$, $w_1\cdot\widetilde{w}_1=(-1)\cdot w_1$ and $w_1\cdot \widetilde{w}_2=0$.\end{center}
\end{prop}

\begin{dem}
By \cite{garibaldi2003}, Theorem 12.3, we just have to check the equalities on the versal $\mathbb D_4$-torsor $T=(K(\sqrt t),u+v\sqrt t)$. Let us prove for instance the first one. We have \begin{equation*}w_2(T)=w_2(\langle 2,2t\rangle)=(2)\cdot (2t)=(2)\cdot (t)= (2)\cdot w_1(\langle 2,2t\rangle).\end{equation*} The other equalities are left to the reader. $\blacksquare$
\end{dem}

\subsection{Cohomological invariants of Weyl groups of type $B_{2m}$}
Let us recall Theorem \ref{icbn} :
\begin{thm*}
Let $k_0$ be a field of characteristic zero, such that $-1$ and $2$ are squares in $k_0$ and let $W$ be a Weyl group of type $B_{2m}$ ($m\geq 1$).  Then the $H^*(k_0,\mathbb Z/2\mathbb  Z)$-module $\textsf{Inv}_{k_0}(W,\mathbb Z/2\mathbb Z)$ is free with basis \begin{equation*}\{w_i\cdot\widetilde{w}_j\}_{0\leq i\leq m,0\leq j\leq 2(m-i)}.\end{equation*}
\end{thm*}

Let us explain the strategy of the proof. We will show that these Stiefel-Whitney invariants generate the module $\textsf{Inv}_{k_0}(W,\mathbb Z/2\mathbb Z)$ by induction on $m\geq 1$. If $m=1$, $W\simeq \mathbb D_4$ and Theorem \ref{d4} gives the answer. Let $m\geq 2$ and let $a\in\textsf{Inv}_{k_0}(W,\mathbb Z/2\mathbb Z)$. Then we will show, as a consequence of Corollary \ref{icbn1} that,  if we consider is a particular subgroup $W_0$ of $W$ isomorphic to $\mathbb D_4\times W'$, where $W'$ is a Weyl group of type $B_{2m-2}$, then $a$ is completely determined by $\textsf{Res}_W^{W_0}(a)$. Since we know the invariants of $\mathbb D_4$ and of $W'$, we can write $\textsf{Res}_W^{W_0}(a)$ in terms of invariants of $W_0$, that we can describe from invariants of $\mathbb D_4$ and $W'$. Then, by a second induction on $0\leq q\leq m$, we will study the restrictions to $H_q$ in order to identify $\textsf{Res}_W^{W_0}(a)$ with restrictions of the required Stiefel-Whitney invariants. In a last part, we will show by induction on $0\leq q\leq m$ that the family is free, computing restrictions of the $w_i\cdot\widetilde{w}_j$ to the subgroups $H_q$.\\

From now on, let $k_0$ be any field of characteristic zero such that $-1$ and $2$ are squares in $k_0$. However, much (but not all) of what follows is true, even if $-1$ and $2$ are no more squares in $k_0$.\\

Let us prove, by induction on $m\geq 1$, that: 
\begin{prop} 
\label{gen}
For any $m\geq 1$, if $W$ is a Weyl group of type $B_{2m}$, the family $\{w_i\cdot \widetilde{w}_j\}_{0\leq i\leq m, 0\leq j\leq 2(m-i)}$ generates $\textsf{Inv}_{k_0}(W,\mathbb Z/2\mathbb Z)$ as an $H^*(k_0,\mathbb Z/2\mathbb Z)$-module.
\end{prop}

If $m=1$, $W\simeq \mathbb D_4$ and Theorem \ref{d4} allows us to conclude.\\

Let $m\geq 2$ and let $W$ be a Weyl group of type $B_{2m}$. Let us assume that any Weyl group $W'$ of type $B_{2(m-1)}$ satisfies the induction hypothesis, i.e. its cohomological invariants form a free $H^*(k_0,\mathbb Z/2\mathbb Z)$-module with a basis given by the invariants $w_i^{W'}\cdot \widetilde{w}_j^{W'}$, with $0\leq i\leq m-1$ and $0\leq j\leq 2(m-1-i)$. \\

Let $W$ be a Weyl group of type $B_{2m}$. It corresponds to a root system \begin{equation*}S=\{\pm e_i,\pm e_i\pm e_j, 1\leq i\leq 2m, 1\leq j\neq i\leq 2m\}\end{equation*}
Let us denote by $W'$ the subgroup of $W$ corresponding to the root subsystem \begin{equation*}S'=\{\pm e_i,\pm e_i\pm e_j, 3\leq i\leq 2m, 3\leq j\neq i\leq 2m\}.\end{equation*} It is a Weyl group of type $B_{2m-2}$. Let us also denote by $W_0$ the non-irreducible Weyl group corresponding to the root subsystem $\{\pm e_1,\pm e_2,\pm e_1\pm e_2\}\sqcup S'$. Then $W_0$ is a subgroup of $W$ isomorphic to  $\mathbb D_4\times W'$.

\begin{lem}
\label{w0}
Any cohomological invariant of $W$ over $k_0$ with coefficients in $\mathbb Z/2\mathbb Z$ is completely determined by its restriction to $W_0$.
\end{lem}

\begin{dem}
By Corollary \ref{icbn1}, any invariant of $W$ is completely determined by its restrictions to the subgroups $H_q$ for $0\leq q\leq m$. Let $0\leq q\leq m$. The root system $S_q$ defined in Section \ref{s1} corresponding to $H_q$ is clearly a subset of $\{\pm e_1,\pm e_2,\pm e_1\pm e_2\}\sqcup S'$. Hence, $H_q\subset W_0$. It implies that, if $k/k_0$ is an extension and if $T_q$ is a $W$-torsor over $k$ which lies in the image of $H^1(k,H_q)\rightarrow H^1(k, W)$, then $T_q$ also lies in the image of $H^1(k,W_0)\rightarrow H^1(k,W)$.$\blacksquare$
\end{dem}

\subsubsection{Restriction of Stiefel-Whitney classes of trace forms and twisted trace forms}
The aim of this section is to give formulae for restrictions of Stiefel-Whitney classes to the subgroups $W_0$ and $H_q$, $q=0,...,m$. Note that $w_0=1=\widetilde{w}_0$. Let us first compute the restrictions of Stiefel-Whitney invariants to $W_0$.
\begin{prop}
\label{siw0}
We have the following formulae :
\begin{itemize}
\item[(i)] for $1\leq j\leq 2m$, $\textsf{Res}_W^{W_0}(\widetilde{w}_j)=\widetilde{w}_2^{\mathbb D_4}\cdot \widetilde{w}_{j-2}^{W'}+ \widetilde{w}_1^{\mathbb D_4}\cdot \widetilde{w}_{j-1}^{W'}+ \widetilde{w}_j^{W'}$;
\item[(ii)] for $1\leq i\leq m$, $\textsf{Res}_W^{W_0}(w_i)=w_i^{W'} + w_1^{\mathbb D_4}\cdot w_{i-1}^{W'}$;
\item[(iii)] for $1\leq i\leq m$ and $1\leq j \leq 2(m-i)$,
\end{itemize}
 \begin{equation*}\begin{aligned}\textsf{Res}_W^{W_0}(w_i\cdot \widetilde{ w}_j)=\widetilde{w}_2^{\mathbb D_4}\cdot w_i^{W'}\cdot\widetilde{w}_{j-2}^{W'} + &\widetilde{w}_1^{\mathbb D_4}\cdot w_i^{W'}\cdot \widetilde{ w}_{j-1}^{W'} \\ &+  w_1^{\mathbb D_4}\cdot w_{i-1}^{W'} \cdot \widetilde{w}_j^{W'}+ w_i^{W'}\cdot \widetilde{ w}_j^{W'}.\end{aligned}\end{equation*}
\end{prop}

\begin{dem}

Let $k/k_0$ be a field extension. Let $(L,\alpha)\in H^1(k,W_0)$. Then $L=L_1\times L_2$ with $L_1$ an \'etale $k$-algebra of rank $2$ and $\alpha=(\alpha_1,\alpha_2)$. Thus, the quadratic form $q_{L,\alpha}:x\mapsto \textsf{Tr}_L(\alpha x^2)$ decomposes into $q_{L,\alpha}=q_{L_1,\alpha_1}\oplus q_{L_2,\alpha_2}$. Hence, for $0\leq j\leq 2m$, \begin{equation*} w_j(q_{L,\alpha})=\underset{0\leq i\leq j}{\sum} w_i(q_{L_1,\alpha_1})\cdot w_{j-i}(q_{L_2,\alpha_2}).\end{equation*} Since $w_i(q_{L_1,\alpha_1})=0$ as soon as $i>2$, we get that \begin{equation*}w_j(q_{L,\alpha})= w_2(q_{L_1,\alpha_1})\cdot w_{j-2}(q_{L_2,\alpha_2})+ w_1(q_{L_1,\alpha_1})\cdot w_{j-1}(q_{L_2,\alpha_2})+  w_j(q_{L_2,\alpha_2})\end{equation*} which gives us \textit{(i)}. 
Likewise, with the quadratic form $q_L:x\mapsto \textsf{Tr}_L(x^2)$, we get that, for $0\leq i\leq m$, \begin{equation*}w_i(q_L)= w_2(q_{L_1})\cdot w_{i-2}(q_{L_2})+ w_1(q_{L_1})\cdot w_{i-1}(q_{L_2})+  w_i(q_{L_2})\end{equation*}  Thus, \textit{(ii)} follows from Proposition \ref{reld4} using the assumption that $-1$ and $2$ are squares in $k_0^\times$ . Since the cup-product commutes with the restriction map, Formula \textit{(iii)} follows from \textit{(i)}, \textit{(ii)} and from Proposition \ref{reld4}. $\blacksquare$
\end{dem}

Let us now consider the restrictions of Stiefel-Whitney invariants to the subgroups $H_q$, for $0\leq q\leq m$.
We do not need here the exhaustive list of all the restrictions, so we only give those that will be useful in the sequel.\\
\begin{lem}
\label{resh_q}
Let $q\in \{0,...,m-1\}$. For all $q+1\leq i \leq m$ and all $0\leq j\leq 2(m-i)$,
\begin{equation*}\textsf{Res}_W^{H_q}(w_i\cdot \widetilde{w}_j)=0.\end{equation*}
\end{lem}

\begin{dem}
 If $k$ is an extension of $k_0$ and if \begin{equation*}T_q=(k(\sqrt{t_1})\times...\times k(\sqrt{t_q})\times k^{2(m-q)}, (u_1,...,u_q,u_{q+1},v_{q+1},...,u_m,v_m))\end{equation*} is a $W$-torsor over $k$ lying in the image of $H^1(k,H_q)\rightarrow H^1(k,W)$, then\\ \[\begin{aligned} w_i(T_q)&=w_i(\langle 2,2t_1,2,2t_2,...,2,2t_q\rangle) \\&=w_i(\langle 1,t_1,...,1,t_q\rangle) \\ &=w_i(\langle t_1,...,t_q\rangle)\hspace{2cm}\end{aligned}\] which is $0$ since $i\geq q+1$. $\blacksquare$
\end{dem}

Let us go further for the case $q=0$. Recall that, for any $I\subset\{1,...,2m\}$, $a_I$ denotes the invariant of $H_0$ given by $(x_1,...,x_{2m})\mapsto (x)_I$, where $(x)_I$ is the cup-product of the $(x_i)$ for $i\in I$ (see Example \ref{ex1}). Recall also that $a_j^{(0)}$ denotes the invariant $\underset{I\subset\{1,...,2m\}; \mid I\mid=j}{\sum} a_I$, for $0\leq j\leq 2m$ (see Example \ref{ex2}).\\
\begin{lem}
\label{h0}
For any $0\leq j\leq 2m$, $\textsf{Res}_W^{H_0}(\widetilde{w}_j)=a_j^{(0)}$. In particular, the family $\{\textsf{Res}_W^{H_0}(\widetilde{w}_j)\}_{0\leq j\leq 2m}$ is free over $H^*(k_0,\mathbb Z/2\mathbb Z)$.\\
\end{lem}

\begin{dem}
If $k/k_0$ is an extension and $T_0=(k^{2m},(u_1,v_1,...,u_m,v_m))$ is a $W$-torsor over $k$ lying in the image of $H^1(k,H_0)\rightarrow H^1(k, W)$, we have \begin{equation*}\widetilde{w}_j(T_0)=w_j(\langle u_1,v_1,...,u_m,v_m\rangle)=a_j^{(0)}(u_1,v_1,...,u_m,v_m).\end{equation*} Moreover, by Proposition \ref{e2}, the invariants $a_j^{(0)}$ form a basis of the submodule $\textsf{Inv}_{k_0}(H_0,\mathbb Z/2\mathbb Z)^{N_0/H_0}$ and this gives the freedom of  $\{\textsf{Res}_W^{H_0}(\widetilde{w}_j)\}_{0\leq j\leq 2m}$. $\blacksquare$
\end{dem}

\begin{lem}
\label{j>2(m-q)}
Let $0\leq i\leq m$.
Then, for any $j>2(m-i)$, $w_i\cdot \widetilde{w}_j=0$.
\end{lem}

\begin{dem}
Let $j>2(m-i)$. By Corollary \ref{icbn1}, it is enough to show that the restriction of $w_i\cdot \widetilde{w}_j$ to any subgroup $H_q$ of $W$ ($0\leq q\leq m$) is zero. Let $0\leq q\leq m$, let $k/k_0$ be a field extension and let $T_q=(k(\sqrt {t_1})\times ...\times k(\sqrt{t_q})\times k^{2(m-q)}, (u_1,...,u_q,u_{q+1},v_{q+1},...,u_m,v_m))$ be a $W$-torsor over $k$ lying in the image of $H^1(k, H_q)\rightarrow H^1(k,W)$. We then have to show that $w_i(T_q)\cdot \widetilde{w}_j(T_q)=0$. By Lemma \ref{resh_q}, $w_i(T_q)=0$ if $q<i$. Let us assume that $q\geq i$. We have :
\begin{equation*}\begin{aligned}w_i(T_q)\cdot \widetilde{w}_j(T_q) &= w_i(\langle t_1,...,t_q\rangle) \cdot  w_j(\langle u_1,u_1t_1,...,u_q,u_qt_q,u_{q+1},v_{q+1}, ...,u_m,v_m\rangle)\\ &=\underset{1\leq j_1<...<j_i \leq q}{\sum} (t_{j_1})\cdot \ldots \cdot (t_{j_i}) \cdot \big[ \sum_{j'=0}^j w_{j'}(\langle u_1,u_1t_1,...,u_q,u_qt_q\rangle)\\ &\hspace{5cm}\cdot w_{j-j'}(\langle u_{q+1},v_{q+1},...,u_m,v_m\rangle) \big].\end{aligned}\end{equation*}

Since the quadratic form $\langle u_{q+1},v_{q+1},...,u_m,v_m\rangle$ has rank $2(m-q)$, we have $w_{j-j'}(\langle u_{q+1},v_{q+1},...,u_{m},v_{m}\rangle)=0$ if $j-j'>2(m-q)$. So we get :
\begin{equation*}\begin{aligned} w_i(T_q)\cdot \widetilde{w}_j(T_q) &=\underset{1\leq j_1<...<j_i \leq q}{\sum} (t_{j_1})\cdot \ldots \cdot (t_{j_i})\\ &\hspace{1cm}\cdot \big[ \sum_{j'=j-2(m-q)}^j w_{j'}(\langle u_1,u_1t_1,...,u_q,u_qt_q\rangle)\\ &\hspace{2cm} \cdot w_{j-j'}(\langle u_{q+1},v_{q+1},...,u_{m},v_{m}\rangle) \big].\end{aligned}\end{equation*}

Since $j>2(m-i)$, then $j-2(m-q)>2(q-i)$, which gives us :\\
\begin{equation}
\label{eq1}
\begin{aligned}w_i(T_q)\cdot \widetilde{w}_j(T_q)&=\underset{1\leq j_1<...<j_i \leq q}{\sum} (t_{j_1})\cdot \ldots \cdot (t_{j_i})\\ &\hspace{2cm}\cdot \big[ \sum_{j'=2(q-i)+1}^j w_{j'}(\langle u_1,u_1t_1,...,u_q,u_qt_q \rangle)\\ &\hspace{3cm} \cdot w_{j-j'}(\langle u_{q+1},v_{q+1},...,u_{m},v_{m}\rangle)\big].\end{aligned}\end{equation}

Let us show that, for any $0\leq j\leq q$ and any  $2(q-i)<j'\leq j$, \begin{equation}\label{eq2}\begin{aligned}(t_j)&\cdot  w_{j'}(\langle u_1,u_1t_1,...,u_q,u_qt_q \rangle)\\ &=(t_j)\cdot w_{j'}(\langle u_1,u_1t_1,...,u_{j-1},u_{j-1}t_{j-1},u_{j+1},u_{j+1}t_{j+1},...,u_q,u_qt_q\rangle).\end{aligned}\end{equation}

\noindent Let  $2(q-i)<j'\leq j$. For sake of simplicity, let us assume that $j=1$. We have
\begin{equation*} \begin{aligned} w_{j'}(\langle u_1,u_1t_1,...,u_q,u_qt_q\rangle)=&(u_1)\cdot (u_1t_1) \cdot w_{j'-2}(\langle u_2,u_2t_2,...,u_q,u_qt_q\rangle)\\ & +(u_1)\cdot  w_{j'-1}(\langle u_2,u_2t_2,...,u_q,u_qt_q\rangle)\\ & + (u_1t_1)\cdot  w_{j'-1}(\langle u_1,u_1t_1,...,u_q,u_qt_q\rangle)\\ & +  w_{j'}(\langle u_2,u_2t_2,...,u_q,u_qt_q\rangle),\end{aligned}\end{equation*} so
\begin{equation*}
\begin{aligned} w_{j'}(\langle u_1,u_1t_1,...,u_q,u_qt_q \rangle) = &(u_1)\cdot (u_1t_1) \cdot w_{j'-2}(\langle u_2,u_2t_2,...,u_q,u_qt_q\rangle)\\ & +(t_1)\cdot  w_{j'-1}(\langle u_2,u_2t_2,...,u_q,u_qt_q\rangle) \\ &+  w_{j'}(\langle u_2,u_2t_2,...,u_q,u_qt_q\rangle).\end{aligned}\end{equation*}

Hence,
\begin{equation*}\begin{aligned} (t_1) \cdot & w_{j'}(\langle u_1,u_1t_1,...,u_q,u_qt_q\rangle)\\ &= (t_1)\cdot (u_1)\cdot (u_1t_1) \cdot w_{j'-2}(\langle u_2,u_2t_2,...,u_q,u_qt_q\rangle ) \\ &\hspace{0.5cm} +(t_1)\cdot (t_1)\cdot  w_{j'-1}(\langle u_2,u_2t_2,...,u_q,u_qt_q\rangle)\\ & \hspace{0.5cm} + (t_1)\cdot w_{j'}(\langle u_2,u_2t_2,...,u_q,u_qt_q\rangle ).\end{aligned}\end{equation*} 

Since $(t_1)\cdot (u_1) \cdot (u_1t_1)=0$ and $(t_1)\cdot (t_1)=(t_1)\cdot(-1)=0$, we get that 
\begin{equation*} (t_1)\cdot w_{j'}(\langle u_1,u_1t_1,...,u_q,u_qt_q\rangle)=(t_1)\cdot w_{j'}(\langle u_2,u_2t_2,...,u_q,u_qt_q\rangle).\end{equation*}
This proves (\ref{eq2}). An immediate induction shows that, for any $0\leq j_1<...<j_i\leq q$, \begin{equation*}\begin{aligned} &(t_{j_1})\cdot\ldots\cdot (t_{j_i})\cdot w_{j'}(\langle u_1,u_1t_1,...,u_q,u_qt_q\rangle)\\ &\hspace{1cm}=(t_{j_1})\cdot\ldots\cdot (t_{j_i})\cdot w_{j'}(\langle u_{j'_1},u_{j'_1}t_{j'_1},...,u_{j'_{q-i}},u_{j'_{q-i}}t_{j'_{q-i}}\rangle)\end{aligned}\end{equation*} where $\{j'_1,...,j'_{q-i}\}$ is the complementary of $\{j_1,...j_i\}$ in $\{1,...,q\}$. Since the quadratic form $Q=\langle u_{j'_1},u_{j'_1}t_{j'_1},...,u_{j'_{q-i}},u_{j'_{q-i}}t_{j'_{q-i}}\rangle$ has rank $2(q-i)$, we get, for any $j'>2(q-i)$, that $w_{j'}(Q)=0$.
Using this in Equation (\ref{eq1}), we can conclude that $w_i(T_q)\cdot \widetilde{w}_j(T_q)=0$. $\blacksquare$
\end{dem}

\begin{rmk}
This lemma does not hold anymore if we do not assume that $-1$ or $2$ are squares in $k_0$.
\end{rmk}

Let us state the last lemma of this section :
\begin{lem}
\label{libHq}
Let $0\leq q\leq m$.
The family $\{\textsf{Res}_W^{H_q}(w_q\cdot \widetilde{w}_j)\}_{0\leq j\leq 2(m-q)}$ is free over $H^*(k_0,\mathbb Z/2\mathbb Z)$.
\end{lem}

\begin{dem}
To show that this family of invariants is free, it is enough to prove it for their value on a versal $H_q$-torsor over $k_0$ (see \cite{garibaldi2003}, Theorem 12.3). Let $t_1,...,t_q,u_1,...,u_m,v_{q+1},...,$ $v_m$ be independent indeterminates over $k_0$ and $K=k_0(t_1,...,t_q,u_1,...,u_m,v_{q+1},...,v_m)$. Let us denote by $T_q$ the image of the versal $H_q$-torsor \begin{equation*}(u_1,u_1t_1,...,u_q,u_qt_q,u_{q+1},v_{q+1},...,u_m,v_m)\end{equation*} by $H^1(K,H_q)\rightarrow H^1(K,W)$. We have to show that the invariants $w_q(T_q)\cdot \widetilde{w}_j(T_q)$ where $0\leq j\leq 2(m-q)$ form a free family over $H^*(k_0,\mathbb Z/2\mathbb Z)$.\\ 

Let $0\leq j\leq 2(m-q)$. We have 
$w_q(T_q)=w_q(\langle t_1,...,t_q\rangle)=(t_1)\cdot \ldots \cdot(t_q)$ and
\begin{equation*}\begin{aligned}\widetilde{w}_j(T_q)&=w_j(\langle 2u_1,2u_1t_1,...,2u_q,2u_qt_q, u_{q+1},v_{q+1},...,u_m, v_m\rangle)\\ &=w_j(\langle u_1,u_1t_1,...,u_q,u_qt_q, u_{q+1},v_{q+1},...,u_m, v_m\rangle)\\ &= \underset{0\leq j'\leq j}{\sum} w_{j'}(\langle u_1,u_1t_1,...,u_q,u_qt_q\rangle)\cdot w_{j-j'}(\langle u_{q+1},v_{q+1},...,u_m,v_m\rangle).\end{aligned}\end{equation*} 
As shown in the proof of Lemma \ref{j>2(m-q)}, for any $1\leq j'\leq j$, \begin{equation*}(t_1)\cdot w_{j'}(\langle u_1,u_1t_1,...,u_q,u_qt_q\rangle)=(t_1)\cdot w_{j'}(\langle u_2,u_2t_2,...,u_q,u_qt_q\rangle).\end{equation*} Thus an easy induction shows that $(t_1)\cdot \ldots \cdot (t_q)\cdot w_{j'}(\langle u_1,u_1t_1,...,u_q,u_qt_q\rangle)=0$.
Hence, \begin{equation*}w_q(T_q)\cdot \widetilde{w}_j(T_q)=(t_1)\cdot \ldots \cdot (t_q) \cdot w_j(\langle u_{q+1},v_{q+1},...,u_m,v_m\rangle).\end{equation*}

Furthermore, since the monomials in $(t_1),...,(t_q),(u_1),...,(u_m),(v_{q+1}),...,(v_m)$ form a free family over $H^*(k_0,\mathbb Z/2\mathbb Z)$, it is easily seen that the invariants $w_q(T_q) \cdot \widetilde{w}_j(T_q)$, for $0\leq j\leq 2(m-q)$, also form a free family over $H^*(k_0,\mathbb Z/2\mathbb Z)$. $\blacksquare$
\end{dem}

\begin{rmk}
This lemma, contrary to Lemma \ref{j>2(m-q)}, is still true if we do not assume anymore that $-1$ or $2$ are squares in $k_0$.
\end{rmk}

\subsubsection{Cohomological invariants of $\textsf{Res}_W^{W_0}(a)$}
Let $a\in\textsf{Inv}_{k_0}(W,\mathbb Z/2\mathbb Z)$. Let us summarize what we got.
By Corollary \ref{icbn1}, the invariant $a$ is completely determined by its values on the $W$-torsors that lie in the image of a map $H^1(k,H_q)\rightarrow H^1(k,W)$ (for $0\leq q \leq m$). In fact, Lemma \ref{w0} yields that $a$ is completely determined by its restriction to the subgroup $W_0$ and so by its values on the $W$-torsors that are the image of a $W_0$-torsor. For any extension $k/k_0$, such a torsor is a $W$-torsor over $k$ of the form $T_1\times T_2$, where $T_1$ is a $\mathbb D_4$-torsor over $k$ and $T_2$ a $W'$-torsor over $k$. In the sequel, we will then work with these $W$-torsors of the form $T_1\times T_2$ over any extension $k/k_0$.\\

Let $k/k_0$ be a field extension and let $T_2$ be a $W'$-torsor over $k$. For any $k'/k$, let us consider the map \begin{equation*} \begin{aligned} (a^{(k,T_2)})_{k'}:H^1(k',\mathbb D_4)&\rightarrow H^*(k',\mathbb Z/2\mathbb Z)\\ T_1&\mapsto a_{k'}(T_1\times \textsf{Res}_{k'/k}(T_2)).\end{aligned}\end{equation*} It is clear that these maps define a cohomological invariant $a^{(k,T_2)}\in \textsf{Inv}_{k}(\mathbb D_4,\mathbb Z/2\mathbb Z)$. By Theorem \ref{d4}, there exist, for $i\in\{0,1,\widetilde 1,\widetilde 2\}$, $(c_i)_k(T_2)\in H^*(k,\mathbb Z/2\mathbb Z)$ such that 
\begin{equation*}a^{(k,T_2)}=(c_0)_k(T_2)+(c_1)_k(T_2)\cdot w_1^{\mathbb D_4} + (c_{\widetilde 1})_k(T_2)\cdot \widetilde{w}_1^{\mathbb D_4}+(c_{\widetilde 2})_k(T_2)\cdot \widetilde{w}_2^{\mathbb D_4}.\end{equation*}

For $i\in\{0,1,\widetilde 1,\widetilde 2\}$ and for $k/k_0$, let us denote by $(c_i)_k:H^1(k,W')\rightarrow H^*(k,\mathbb Z/2\mathbb Z)$ the induced map.

\begin{prop}
\label{ci}
If $i\in\{0,1,\widetilde 1,\widetilde 2\}$, $c_i$ is a cohomological invariant of $W'$ over $k_0$.
\end{prop}

The proof of Proposition \ref{ci} is left to the reader. The induction hypothesis of Proposition \ref{gen} implies that, for every $l\in\{0,1,\widetilde 1,\widetilde 2\}$, there exists $b_{l,i,j}\in H^*(k_0,\mathbb Z/2\mathbb Z)$, for any $0\leq i\leq m-1$ and  any $0\leq j\leq 2(m-1-i)$ such that \begin{equation*}c_l=\underset{0\leq i\leq m-1, 0\leq j\leq 2(m-1-i)}{\large\sum} b_{l,i,j}\cdot w_i^{W'}\cdot \widetilde{w}_j^{W'}.\end{equation*}

From now on, we will use sometimes the notation $w_{\widetilde 1}^{\mathbb D_4}=\widetilde{w}_1^{\mathbb D_4}$ and $w_{\widetilde 2}^{\mathbb D_4}=\widetilde{ w}_2^{\mathbb D_4}$ in order to simplify the sums to write.
Hence we have :
\begin{equation}\textsf{Res}_W^{W_0}(a)= \underset{l\in\{0,1, \widetilde 1,\widetilde 2\}, 0\leq i\leq m-1, 0\leq j\leq 2(m-1-i)}{\large\sum} b_{l,i,j}\cdot w_l^{\mathbb D_4}\cdot w_i^{W'}\cdot \widetilde{w}_j^{W'}.\end{equation}

\subsubsection{Restrictions of $\textsf{Res}_W^{W_0}(a)$ to $H_q$, for $0\leq q \leq m$}
We will now show the following proposition by induction on $q\in\{0,...,m-1\}$:

\begin{prop}
\label{rec}
There are some coefficients $C_{i,j}\in H^*(k_0,\mathbb Z/2\mathbb Z)$ such that, for any $0\leq q\leq m-1$, 
\begin{equation*}\textsf{Res}_W^{W_0}(a)=\underset{0\leq i\leq q, 0\leq j\leq 2(m-i)}{\sum} C_{i,j}\cdot \textsf{Res}_W^{W_0}(w_i\cdot \widetilde{w}_j) +
a_{q+1}\end{equation*} where \begin{equation*} \begin{aligned} a_{q+1}=&\underset{l\in\{0,\widetilde 1,\widetilde 2\},q+1\leq i\leq m-1, 0\leq j \leq 2(m-1-i)}{\sum} b_{l,i,j}\cdot w_l^{\mathbb D_4}\cdot w_i^{W'}\cdot \widetilde{w}_j^{W'}\\ &+ \underset{q\leq i <m-1, 0\leq j\leq 2(m-1-i)} {\sum} b_{1,i,j}\cdot w_1^{\mathbb D_4}\cdot w_i^{W'}\cdot \widetilde{w}_j^{W'}.\end{aligned}\end{equation*}\\
\end{prop}

In other words, at each step $q$ of our induction, we will identify parts of the sums with a linear combination of the invariants $w_q\cdot \widetilde{w}_j$ for $0\leq j\leq 2(m-q)$ by considering the restriction to the subgroup $H_q$, where we have a lot of information about torsors, Stiefel-Whitney invariants, etc. This will then reduce the extra term $a_{q+1}$. We will finally show that at rank $m$ of the induction, the extra term will have completely disappeared.

\begin{dem}
Let us check the case $q=0$ first : let us consider the restriction $\textsf{Res}_W^{H_0}(a)$. Let $k/k_0$ be an extension, let $T_1=(k^2,(u_1,v_1))$ be a $\mathbb D_4$-torsor over $k$ and let $T'=(k^{2m-2},(u_2,v_2,...,u_m,v_m))$ be a $W'$-torsor over $k$ so that the cohomology class associated with $T_1\times T'$ lies in the image of $H^1(k,H_0)\rightarrow H^1(k,W)$.  Hence,
\begin{equation*}a_k(T_1\times T')=\underset{ 0\leq i\leq m-1, 0\leq j\leq 2(m-1-i)}{\underset{l\in\{0,1,\widetilde 1,\widetilde 2\},}{\large\sum}} b_{l,i,j}\cdot w_l^{\mathbb D_4}(T_1)\cdot w_i^{W'}(T')\cdot \widetilde{ w}_j^{W'}(T').\end{equation*}

For any $1\leq i\leq m-1$, $w_i^{W'}(T')=0$ and $w_1^{\mathbb D_4}(T_1)=0$, so we have :
\begin{equation*}a_k(T_1\times T')=\underset{l\in \{0,\widetilde 1,\widetilde 2\}, 0\leq j\leq 2(m-1)}{\large\sum} b_{l,0,j}\cdot w_l^{\mathbb D_4}(T_1)\cdot \widetilde{w}_j^{W'}(T').\end{equation*}

Let us embed $H_0$ in $W_0=\mathbb D_4\times W'$. Then $H_0$ decomposes in this product into two factors, the left one being isomorphic to $(\mathbb Z/2\mathbb Z)^2$ and denoted by $H_0^{\mathbb D_4}$ and the right factor being an abelian subgroup of $W'$ generated by reflections. We denote it by $H'_0$. Note that $H'_0$ is for $W'$ exactly what $H_0$ is for $W$.  Lemma \ref{h0} applies here for $W'$ and $H_0'$ : \begin{equation*}\textsf{Res}_{W'}^{H_0'}(\widetilde{w}_j^{W'})=a_j^{(0)}=\underset{J\subset \{3,...,2m\}, \mid J\mid=j}{\sum}a_J.\end{equation*}

For any $0\leq j\leq 2(m-1)$, we have \begin{equation*}\begin{aligned} 
\textsf{Res}_{\mathbb D_4}^{H_0^{\mathbb D_4}}(\widetilde{w}_1^{\mathbb D_4})\cdot  \textsf{Res}_{W'}^{H_0'}(\widetilde{w}_j^{W'}) &=\big( a_{\{1\}}+a_{\{2\}} \big)\cdot \big( \underset{J\subset \{3,...,2m\}, \mid J \mid=j}{\sum} a_J \big)\\ &= \underset{J\subset \{3,...,2m\}, \mid J \mid=j}{\sum} \big(a_{ \{1\} \cdot J} + a_{\{2\}\cdot J} \big)\end{aligned}\end{equation*}  and 
\begin{equation*}\begin{aligned}\textsf{Res}_{\mathbb D_4}^{H_0^{\mathbb D_4}}(\widetilde{w}_2^{\mathbb D_4})\cdot  \textsf{Res}_{W'}^{H_0'}(\widetilde{w}_j^{W'}) &=a_{\{1,2\}}\cdot \big(\underset{J\subset \{3,...,2m\}, \mid J \mid=j}{\sum} a_J\big)\\ &=\underset{J\subset \{3,...,2m\}, \mid J \mid=j}{\sum} a_{\{1,2\}\cdot J}.\end{aligned}\end{equation*}

Let us come back to $\textsf{Res}_W^{H_0}(a)$. We get that
\begin{equation*}\begin{aligned}\textsf{Res}_W^{H_0}(a) =\underset{0\leq j\leq 2(m-1)}{\large\sum} \big[ & b_{0,0,j}\cdot \big(\underset{J\subset \{3,...,2m\}, \mid J \mid=j}{\sum} a_J\big)\\ +  &b_{\widetilde 1,0,j}\cdot \big(\underset{J\subset \{3,...,2m\}, \mid J \mid=j}{\sum} (a_{J\cdot \{1\}} + a_{J \cdot \{2\}})\big)\\  + &b_{\widetilde 2,0,j}\cdot \big(\underset{J\subset \{3,...,2m\}, \mid J \mid=j}{\sum} a_{J\cdot \{1,2\}}\big)\big] .\end{aligned}\end{equation*}

Moreover, $\textsf{Res}_W^{H_0}(a)$ belongs to the submodule of the cohomological invariants of $H_0$ fixed by the group $N_0/H_0$. By Proposition \ref{e2}, for any $i=0,...,2m$, there exists $b_i\in H^*(k_0,\mathbb Z/2\mathbb Z)$ such that \begin{equation*}\textsf{Res}_W^{H_0}(a)=\overset{2m}{\underset{i=0}{\sum}} b_i\cdot a_i^{(0)},\end{equation*} where, for $i=0,...,2m$, $a_i^{(0)}=\underset {I\subset \{1,...,2m\}, \mid I\mid=i}{\sum} a_I$.\\ 

Furthermore, the family $(a_I)_{I\subset \{1,...,2m\}}$ is free in $\textsf{Inv}_{k_0}(H_0,\mathbb Z/2\mathbb Z)$ (see Proposition \ref{e1}), so we get the following relations : for any $0\leq j\leq 2(m-1)$, \begin{center}
$b_{0,0,j}=b_j$, $b_{\widetilde 1,0,j}=b_{j+1}$ and $b_{\widetilde 2,0,j}=b_{j+2}$.\end{center} We can now say that \begin{equation} \label{eq3} \begin{aligned} b_{0,0,1}&=b_{\widetilde 1,0,0},\\ \textrm{ for any } j\geq 2, b_{0,0,j}&=b_{\widetilde 1,0,j-1}=b_{\widetilde 2,0,j-2} \textrm{ and }\\  b_{\widetilde 1,0,2m-2}&=b_{\widetilde 2,0,2m-3}.\end{aligned}\end{equation}

If we set $a'_0=\textsf{Res}_W^{W_0}(a)+a_1$ (it is a cohomological invariant of $W_0$), then 
\begin{equation*} \begin{aligned} a'_0&=\underset{l\in\{0,\widetilde 1,\widetilde 2\}, 0\leq j \leq 2(m-1)}
{\sum} b_{l,0,j}\cdot w_l^{\mathbb D_4}\cdot \widetilde{w}_j^{W'}\\  &=b_{0,0,0}+ b_{0,0,1}\cdot \widetilde{w}_1^{W'} + 
b_{\widetilde 1, 0,0}\cdot \widetilde{w}_1^{\mathbb D_4}\\ & \hspace{0.5cm} +\sum_{j=2}^{2m-2} \big( b_{0,0,j}\cdot \widetilde{w}_j^{W'}+ b_{\widetilde 1,0,j-1}\cdot\widetilde{w}_1^{\mathbb D_4}\cdot \widetilde{w}_{j-1}^{W'} + b_{\widetilde 2,0,j}\cdot\widetilde{w}_2^{\mathbb D_4}\cdot \widetilde{w}_{j-2}^{W'}\big)\\ &\hspace{0.5cm} + b_{\widetilde 1,0,2m-2}\cdot \widetilde{w}_1^{\mathbb 
D_4}\cdot \widetilde{w}_{2m-2}^{W'} + b_{\widetilde 2,0,2m-3}\cdot \widetilde{w}_2^{\mathbb D_4}\cdot \widetilde{w}_{2m-3}^{W'} \\ 
&\hspace{0.5cm}+b_{\widetilde 2,0,2m-2}\cdot \widetilde{w}_2^{\mathbb D_4}\cdot \widetilde{w}_{2m-2}^{W'}.\end{aligned}
\end{equation*}

Therefore, using Relations (\ref{eq3}), we get :
\begin{equation*}
\begin{aligned}
a'_0=&b_{0,0,0}+ b_{0,0,1}\cdot \big(\widetilde{w}_1^{W'} + \widetilde{w}_1^{\mathbb D_4}\big)\\ & +\sum_{j=2}^{2m-2} 
b_{0,0,j}\cdot \big( \widetilde{w}_j^{W'}+ \widetilde{w}_1^{\mathbb D_4}\cdot \widetilde{w}_{j-1}^{W'} + \widetilde{w}_2^{\mathbb 
D_4}\cdot \widetilde{w}_{j-2}^{W'}\big) \\ &+ b_{\widetilde 1,0,2m-2}\cdot \big(\widetilde{w}_1^{\mathbb D_4}\cdot \widetilde{w}_{2m-2}^{W'} + \widetilde{w}_2^{\mathbb D_4}\cdot \widetilde{w}_{2m-3}^{W'} \big)\\ &+ b_{\widetilde 
2,0,2m-2}\cdot \widetilde{w}_2^{\mathbb D_4}\cdot \widetilde{w}_{2m-2}^{W'}.
\end{aligned}
\end{equation*}

By Lemma \ref{siw0}, for $0\leq j\leq 2m$, \begin{equation*}\textsf{Res}_W^{W_0}(\widetilde{w}_j)=\widetilde{w}_2^{\mathbb D_4}\cdot \widetilde{w}_{j-2}^{W'}+ \widetilde{w}_1^{\mathbb D_4}\cdot \widetilde{w}_{j-1}^{W'}+ \widetilde{w}_j^{W'},\end{equation*} thus :
\begin{equation*}\begin{aligned} a'_0=& b_{0,0,0}\cdot \textsf{Res}_W^{W_0}(\widetilde{w}_0)+ b_{0,0,1}\cdot \textsf{Res}_W^{W_0}(\widetilde{w}_1) +\sum_{j=2}^{2m-2} b_{0,0,j}\cdot \textsf{Res}_W^{W_0}(\widetilde{w}_j)\\ & + b_{\widetilde 1,0,2m-2}\cdot \textsf{Res}_W^{W_0}(\widetilde{w}_{2m-1}) + b_{\widetilde 2,0,2m-2}\cdot \textsf{Res}_W^{W_0}(\widetilde{w}_{2m})\\ =& \sum_{j=0}^{2m-2} b_{0,0,j}\cdot \textsf{Res}_W^{W_0}(\widetilde{w}_j) + b_{\widetilde 1,0,2m-2}\cdot \textsf{Res}_W^{W_0}(\widetilde{w}_{2m-1})\\ & + b_{\widetilde 2,0,2m-2}\cdot \textsf{Res}_W^{W_0}(\widetilde{w}_{2m}).
\end{aligned}\end{equation*}

This concludes the case $q=0$.\\

Assume now that $1\leq q \leq m-1$ and that the induction hypothesis is true for the rank $q-1$. By induction hypothesis (see Proposition \ref{rec}), we want to study the extra term $a_q$. Note that $a_q$ is a cohomological invariant of $W_0$. Recall that \begin{equation*}
\begin{aligned} a_q&=\underset{l\in\{0,\widetilde 1,\widetilde 2\},q\leq i\leq m-1, 0\leq j \leq 2(m-1-i)}{\sum} b_{l,i,j}\cdot w_l^{\mathbb D_4}\cdot w_i^{W'}\cdot \widetilde{w}_j^{W'}\\  &+ \underset{q-1\leq i \leq m-1, 0\leq j\leq 2(m-1-i)} {\sum} b_{1,i,j}\cdot w_1^{\mathbb D_4}\cdot w_i^{W'}\cdot \widetilde{w}_j^{W'}. \end{aligned}\end{equation*} 

We then have to show that \begin{equation*}a_q=\underset{0\leq j\leq 2(m-q)}{\sum} C_{q,j}\cdot \textsf{Res}_W^{W_0}(w_q\cdot \widetilde{w}_j) + a_{q+1},\end{equation*} where $C_{q,j}\in H^*(k_0,\mathbb Z/2\mathbb Z)$ for $0\leq j\leq 2(m-q)$. Let $k/k_0$ be an extension, let $T_1=(k^2,(u_1,v_1))$ be a $\mathbb D_4$-torsor over $k$, let $T_2=(k(\sqrt t_2),u_2))$ and let \begin{equation*} T_3=(k(\sqrt {t_3}) \times ...\times k(\sqrt t_{q+1})\times k^{2(m-q-1)}, (u_3,...,u_{q+1},u_{q+2},v_{q+2},...,u_m,v_m))\end{equation*} so that $T_2\times T_3$ is a $W'$-torsor over $k$. Then $T_1\times (T_2\times T_3)$ is a $W_0$-torsor which lies in the image of $H^1(k,H_q)\rightarrow H^1(k, W)$. Since $w_i^{W'}(T_2\times T_3)=w_i(\langle t_2,...,t_{q+1}\rangle)$, we get that $w_i^{W'}(T_2\times T_3)=0$ if $i\geq q+1$. On the other hand, $w_1^{\mathbb D_4}(T_1)=0$. Therefore, we obtain that
\begin{equation}\label{eq22}(a_q)_k(T_1\times (T_2\times T_3))= \underset{0\leq j \leq 2(m-1-q)}{\underset{l\in\{0,\widetilde 1,\widetilde 2\}, }{\sum} }b_{l,q,j}\cdot w_l^{\mathbb D_4}(T_1)\cdot w_q^{W'}(T_2\times T_3)\cdot \widetilde{w}_j^{W'}(T_2\times T_3).\end{equation}

Let us now consider the $\mathbb D_4$-torsor $T_2=(k(\sqrt t_2), u_2)$ and the $W'$-torsor \begin{equation*}\begin{aligned}T_1\times T_3=\big(k^2\times k(\sqrt t_3)\times ...\times & k(\sqrt t_{q+1})\times k^{2(m-q-1)},\\ &(u_1,v_1,u_3,...,u_{q+1},u_{q+2},v_{q+2},...,u_m,v_m)\big).\end{aligned}\end{equation*} 

Then $w_i^{W'}(T_1\times T_3)=w_i(\langle t_3,...,t_{q+1}\rangle)$, so if $i\geq q$, $w_i^{W'}(T_1\times T_3)=0$. Hence,
\begin{equation}\label{eq23}(a_q)_k(T_2\times (T_1\times T_3))=  \underset{0\leq j\leq 2(m-q)} {\sum} b_{1,q-1,j}\cdot w_1^{\mathbb D_4}(T_2)\cdot w_{q-1}^{W'}(T_1\times T_3)\cdot \widetilde{w}_j^{W'}(T_1\times T_3).\end{equation}

Since the two $W_0$-torsors $T_1\times (T_2\times T_3)$ and $T_2\times (T_1\times T_3)$ are isomorphic, it follows from (\ref{eq22}) and (\ref{eq23}) that
\begin{equation} \label{eq4}\begin{aligned} \underset{l\in\{0,\widetilde 1,\widetilde 2\} 0\leq j \leq 2(m-1-q)}{\sum} &b_{l,q,j} \cdot w_l^{\mathbb D_4}(T_1)\cdot w_q^{W'}(T_2\times T_3)\cdot \widetilde{w}_j^{W'}(T_2\times T_3)\\ =\underset{0\leq j\leq 2(m-q)} {\sum} &b_{1,q-1,j}\cdot w_1^{\mathbb D_4}(T_2)\cdot w_{q-1}^{W'}(T_1\times T_3)\cdot \widetilde{w}_j^{W'}(T_1\times T_3).\end{aligned}\end{equation}

Now set  $k_1=k_0(u_2,...,u_m,t_2,...,t_{q+1},v_{q+2},...,v_m)$ and assume that $u_1$ and $v_1$ are independent indeterminates over $k_1$. Then the family $\{1,\widetilde{w}_1^{\mathbb D_4}(T_1),\widetilde{w}_2^{\mathbb D_4}(T_1)\}$ is free over $H^*(k_1,\mathbb Z/2\mathbb Z)$. We then have to collect classes $w_l^{\mathbb D_4}(T_1)$ in (\ref{eq4}). Denoting in an analogous way to $W'\subset W$, by $W''$ the ``same" subgroup of $W'$, we get, by Proposition \ref{siw0} and since  $w_1^{\mathbb D_4}(T_1)=0$, that, for any $0\leq j\leq 2(m-q)$ :  
\begin{equation*}\begin{aligned} w_{q-1}^{W'}(T_1\times T_3)\cdot \widetilde{w}_j^{W'}(T_1\times T_3) = &\widetilde{w}_2^{\mathbb D_4}(T_1)\cdot w_{q-1}^{W''}(T_3)\cdot\widetilde{w}_{j-2}^{W''}(T_3)\\ &+ \widetilde{w}_1^{\mathbb D_4}(T_1)\cdot w_{q-1}^{W''}(T_3)\cdot \widetilde{w}_{j-1}^{W''}(T_3)\\ &+ w_{q-1}^{W''}(T_3)\cdot \widetilde{w}_j^{W''}(T_3).\end{aligned}\end{equation*}

Since the family $\{1,\widetilde{w}_1^{\mathbb D_4}(T_1),\widetilde{w}_2^{\mathbb D_4}(T_1)\}$ is free, we obtain from (\ref{eq4}) the following equalities :
for any $l\in\{0,\widetilde 1,\widetilde 2\}$,
\begin{equation}\label{eq5}\begin{aligned}&\underset{ 0\leq j \leq 2(m-1-q)}{\sum} b_{l,q,j} \cdot w_q^{W'}(T_2\times T_3)\cdot \widetilde{w}_j^{W'}(T_2\times T_3)\\ &=  \underset{0\leq j\leq 2(m-q)} {\sum} b_{1,q-1,j}\cdot w_1^{\mathbb D_4}(T_2)\cdot w_{q-1}^{W''}(T_3)\cdot \widetilde{w}_{j-l}^{W''}(T_3).\end{aligned}\end{equation}

Now set $k_2=k_0(u_3,...,u_m, t_3,...,t_{q+1}, v_{q+2},...,v_m)$ and let $t_2,u_2$ be independent indeterminates over $k_2$. Then the family $\{1,\widetilde{w}_1^{\mathbb D_4}(T_2),\widetilde{w}_2^{\mathbb D_4}(T_2)\}$ is free in $H^*(k_2,\mathbb Z/2\mathbb Z)$ (and $w_1^{\mathbb D_4}(T_2)=\widetilde{w}_1^{\mathbb D_4}(T_2)$). We then have to collect these terms : since $w_q^{W'}(T_3)=0$, by Proposition \ref{siw0}, for any $0\leq j\leq 2(m-1-q)$, we have
\begin{equation*}w_q^{W'}(T_2\times T_3)\cdot \widetilde{w}_j^{W'}(T_2\times T_3)= w_1^{\mathbb D_4}(T_2)\cdot w_{q-1}^{W''}(T_3)\cdot \widetilde{w}_j^{W''}(T_3).\end{equation*}

Hence, for any $l\in\{0,\widetilde 1,\widetilde 2\}$, we get from (\ref{eq5}) that
\begin{equation*}\label{eq6} \begin{aligned}\underset{ 0\leq j \leq 2(m-1-q)}{\sum} b_{l,q,j} &\cdot  w_{q-1}^{W''}(T_3)\cdot \widetilde{w}_j^{W''}(T_3)\\ &=  \underset{0\leq j\leq 2(m-q)} {\sum} b_{1,q-1,j}\cdot w_{q-1}^{W''}(T_3)\cdot \widetilde{w}_{j-l}^{W''}(T_3). \end{aligned}\end{equation*}

We have the three following equalities:\\
for $l=0$ :
\begin{equation}\label{eq7}\begin{aligned}0= &\underset{ 0\leq j \leq 2(m-1-q)}{\sum} \big( b_{0,q,j} + b_{1,q-1,j}\big) \cdot  w_{q-1}^{W''}(T_3)\cdot \widetilde{w}_j^{W''}(T_3)\\ &+ b_{0,q-1, 2(m-q)-1}\cdot w_{q-1}^{W''}(T_3)\cdot \widetilde{w}_{2(m-q)-1}^{W''}(T_3)\\ &+ b_{0,q-1, 2(m-q)}\cdot w_{q-1}^{W''}(T_3)\cdot \widetilde{w}_{2(m-q)}^{W''}(T_3); \end{aligned}\end{equation}

for $l=\widetilde 1$ :
\begin{equation}\label{eq8}\begin{aligned}0=&\underset{ 0\leq j \leq 2(m-1-q)}{\sum} \big( b_{\widetilde 1,q,j} + b_{1,q-1,j+1}\big) \cdot  w_{q-1}^{W''}(T_3)\cdot \widetilde{w}_j^{W''}(T_3)\\ &+ b_{1,q-1,2(m-q)}\cdot w_{q-1}^{W''}(T_3)\cdot \widetilde{w}_{2(m-q)-1}^{W''}(T_3); \end{aligned}\end{equation}

for $l=\widetilde 2$ :
\begin{equation}\label{eq9}0=\underset{ 0\leq j \leq 2(m-1-q)}{\sum} \big( b_{\widetilde 2,q,j} + b_{1,q-1,j+2}\big) \cdot  w_{q-1}^{W''}(T_3)\cdot \widetilde{w}_j^{W''}(T_3).\end{equation}

We now apply Lemma \ref{j>2(m-q)} replacing $W$ by $W''$ and $q$ by $q-1$: \begin{center}$w_{q-1}^{W''}\cdot\widetilde{w}_j^{W''}=0$ if $j>2(m-2-(q-1))=2(m-q)-2$.\end{center} Therefore, \begin{equation*}w_{q-1}^{W''}(T_3)\cdot\widetilde{w}_{2(m-q)-1}^{W''}(T_3)=0=w_{q-1}^{W''}(T_3)\cdot\widetilde{w}_{2(m-q)}^{W''}(T_3).\end{equation*} Thus Relations (\ref{eq7}), (\ref{eq8}) and (\ref{eq9}) become :\\
for $l=0$ :
 \begin{equation}\label{eq10}\underset{ 0\leq j \leq 2(m-1-q)}{\sum} \big( b_{0,q,j} + b_{1,q-1,j}\big) \cdot  w_{q-1}^{W''}(T_3)\cdot \widetilde{w}_j^{W''}(T_3) =0;\end{equation}
for $l=\widetilde 1$ :
 \begin{equation}\label{eq11}\underset{ 0\leq j \leq 2(m-1-q)}{\sum} \big( b_{\widetilde 1,q,j} + b_{1,q-1,j+1}\big) \cdot  w_{q-1}^{W''}(T_3)\cdot \widetilde{w}_j^{W''}(T_3) =0; \end{equation} 

for $l=\widetilde 2$ :
 \begin{equation}\label{eq12} \underset{ 0\leq j \leq 2(m-1-q)}{\sum} \big( b_{\widetilde 2,q,j} + b_{1,q-1,j+2}\big) \cdot  w_{q-1}^{W''}(T_3)\cdot \widetilde{w}_j^{W''}(T_3) =0.\end{equation}

Assume now that $t_3,...,t_{q+1},u_3,...,u_{m},v_{q+2},...,v_m$ are independent indeterminates over $k_0$. Let us denote by $H_q''$ the subgroup of $W''$, defined similarly to $H_q\subset W$. Then, replacing $W$ by $W''$ and $q$ by $q-1$, Lemma \ref{libHq} implies that the invariants $\textsf{Res}_{W''}^{H_q''}(w_{q-1}^{W''}\cdot \widetilde{w}_j^{W''})$ (with $0\leq j\leq 2(m-1-q)$) form a free family over $H^*(k_0,\mathbb Z/2\mathbb Z)$. Since $T_3$ is clearly the image of a versal torsor of $H_q''$, we get, as a direct consequence of \cite{garibaldi2003}, Theorem 12.3, that, for $0\leq j\leq 2(m-1-q)$, the invariants $w_{q-1}^{W''}(T_3)\cdot \widetilde{w}_j^{W''}(T_3)$ form a free family over $H^*(k_0,\mathbb Z/2\mathbb Z)$. It then follows from (\ref{eq10}), (\ref{eq11}) and (\ref{eq12}) that, for any $l\in\{0,1,2\}$ and any $j\in \{0,...,2(m-1-q)\}$, 
\begin{equation}
\label{rel1}
b_{\widetilde l,q,j} + b_{1,q-1,j+l}=0.
\end{equation}

Reordering equalities (\ref{rel1}), we get that, for any $2\leq j\leq 2(m-1-q)$ :
\begin{equation} \label{eq13} \begin{aligned} b_{0,q,0} =& b_{1,q-1,0}, \\
b_{0,q,1}= &b_{\widetilde 1 ,q,0} = b_{1,q-1,1}, \\
b_{\widetilde 2,q,j-2}= b_{\widetilde 1,q,j-1}=& b_{0,q,j} = b_{1,q-1,j},\\
b_{\widetilde 2,q,2(m-1-q)-1}=&b_{\widetilde 1,q,2(m-1-q)}= b_{1,q-1,2(m-q)-1},\\
b_{\widetilde 2,q,2(m-1-q)}=& b_{1,q-1,2(m-q)-1}.\end{aligned}\end{equation}

Let us now come back to $a_q$ :
\begin{equation*}\begin{aligned} a_q &= \underset{l\in\{0,\widetilde 1,\widetilde 2\},q\leq i\leq m-1, 0\leq j \leq 2(m-1-i)}{\sum} b_{l,i,j}\cdot w_l^{\mathbb D_4}\cdot w_i^{W'}\cdot \widetilde{w}_j^{W'}\\ &\hspace{0.5cm} + \underset{q-1\leq i \leq m-1, 0\leq j\leq 2(m-1-i)} {\sum} b_{1,i,j}\cdot w_1^{\mathbb D_4}\cdot w_i^{W'}\cdot \widetilde{w}_j^{W'}, \end{aligned}\end{equation*} so :
\begin{equation*}\begin{aligned} a_q &= \underset{l\in\{0,\widetilde 1,\widetilde 2\}, 0\leq j \leq 2(m-1-q)}{\sum} b_{l,q,j}\cdot w_l^{\mathbb D_4}\cdot w_q^{W'}\cdot \widetilde{w}_j^{W'}\\ & \hspace{0.5cm}+ \underset{0\leq j\leq 2(m-q)} {\sum} b_{1,q-1,j}\cdot w_1^{\mathbb D_4}\cdot w_{q-1}^{W'}\cdot \widetilde{w}_j^{W'}\hspace{1cm}  + a_{q+1}.\end{aligned}\end{equation*}

Set $a'_q=a_q+a_{q+1}$. Using relations (\ref{eq13}), we have  :
\begin{equation}\label{eq14}\begin{aligned} a'_q &=b_{1,q-1,0} \cdot \big( w_q^{W'}+w_1^{\mathbb D_4}\cdot w_{q-1}^{W'} \big)\\ 
&\hspace{0.5cm}+ b_{1,q-1,1}\cdot \big( \widetilde{w}_1^{\mathbb D_4}\cdot w_q^{W'} + w_1^{\mathbb D_4}\cdot w_{q-1}^{W'}\cdot \widetilde{w}_1^{W'} + w_q^{W'}\cdot \widetilde{w}_1^{W'} \big)\\ 
&\hspace{0.5cm}+  \underset{2\leq j \leq 2(m-1-q)}{\sum} b_{1,q-1,j}\cdot \big( \widetilde{w}_2^{\mathbb D_4}\cdot w_q^{W'}\cdot \widetilde{w}_{j-2}^{W'}+ \widetilde{w}_1^{\mathbb D_4}\cdot w_q^{W'}\cdot \widetilde{w}_{j-1}^{W'}\\  &\hspace{5cm} + w_1^{\mathbb D_4}\cdot w_{q-1}^{W'}\cdot \widetilde{w}_j^{W'} + w_q^{W'}\cdot \widetilde{w}_j^{W'}  \big)\\ 
&\hspace{0.5cm}+ b_{1,q-1,2(m-q)-1}\cdot \big(  \widetilde{w}_2^{\mathbb D_4}\cdot w_q^{W'}\cdot \widetilde{w}_{2(m-1-q)-1}^{W'}+ \\ &\hspace{3.5cm}\widetilde{w}_1^{\mathbb D_4}\cdot w_q^{W'}\cdot \widetilde{w}_{2(m-1-q)}^{W'}+ w_1^{\mathbb D_4}\cdot w_{q-1}^{W'}\cdot \widetilde{w}_{2(m-q)-1}^{W'}  \big)\\ 
&\hspace{0.5cm}+ b_{1,q-1,2(m-q)}\cdot \big( \widetilde{w}_2^{\mathbb D_4}\cdot w_q^{W'}\cdot \widetilde{w}_{2(m-1-q)}^{W'} + w_1^{\mathbb D_4}\cdot w_{q-1}^{W'}\cdot \widetilde{w}_{2(m-q)}^{W'} \big). \end{aligned}\end{equation}

Recall now the formulae of Proposition \ref{siw0}. For any $2\leq j\leq 2(m-1-q)$, 
\begin{equation}
\label{eq15}
w_q=  w_q^{W'}+ w_1^{\mathbb D_4}\cdot w_{q-1}^{W'},
\end{equation}

\begin{equation}
\label{eq16}
w_q\cdot \widetilde{w}_1=\widetilde{w}_1^{\mathbb D_4}\cdot w_q^{W'} +  w_1^{\mathbb D_4}\cdot w_{q-1}^{W'} \cdot \widetilde{w}_1^{W'} + w_q^{W'}\cdot \widetilde{w}_1^{W'},\end{equation}

\begin{equation}
\label{eq17}
\begin{aligned} w_q\cdot \widetilde{w}_j = &\widetilde{w}_2^{\mathbb D_4}\cdot w_q^{W'}\cdot\widetilde{ w}_{j-2}^{W'} + \widetilde{w}_1^{\mathbb D_4}\cdot w_q^{W'}\cdot \widetilde{w}_{j-1}^{W'} +  w_1^{\mathbb D_4}\cdot w_{q-1}^{W'} \cdot \widetilde{w}_j^{W'}\\ & + w_q^{W'}\cdot \widetilde{w}_j^{W'},\end{aligned}\end{equation}

\begin{equation}
\label{eq18}
\begin{aligned} w_q\cdot \widetilde{w}_{2(m-q)-1} = &\widetilde{w}_2^{\mathbb D_4}\cdot w_q^{W'}\cdot\widetilde{ w}_{2(m-1-q)-1}^{W'} + \widetilde{w}_1^{\mathbb D_4}\cdot w_q^{W'}\cdot \widetilde{w}_{2(m-1-q)}^{W'}\\  & +  w_1^{\mathbb D_4}\cdot w_{q-1}^{W'} \cdot \widetilde{w}_{2(m-q)}^{W'}\end{aligned}\end{equation}

and
\begin{equation}
\label{eq19}
w_q\cdot \widetilde{w}_{2(m-q)}=\widetilde{w}_2^{\mathbb D_4}\cdot w_q^{W'}\cdot\widetilde{w}_{2(m-1-q)}^{W'} +  w_1^{\mathbb D_4}\cdot w_{q-1}^{W'} \cdot \widetilde{w}_j^{W'}. 
\end{equation}

Note that we used Lemma  \ref{j>2(m-q)} : $w_q^{W'} \cdot \widetilde{w}_{2(m-q)-1}^{W'}=0$ and $w_q^{W'} \cdot \widetilde{w}_{2(m-q)}^{W'}=0$. Therefore, using relations (\ref{eq15}) to (\ref{eq19}) in relation (\ref{eq14}), we get that
\begin{equation}
\begin{aligned} a'_q &=b_{1,q-1,0} \cdot \textsf{Res}_W^{W_0}(w_q) + b_{1,q-1,1}\cdot \textsf{Res}_W^{W_0}(w_q\cdot\widetilde{w}_1)\\ &\hspace{1cm}+  \underset{2\leq j \leq 2(m-1-q)}{\sum} b_{1,q-1,j}\cdot \textsf{Res}_W^{W_0} (w_q\cdot\widetilde{w}_j)\\ &\hspace{1cm} + b_{1,q-1,2(m-q)-1}\cdot \textsf{Res}_W^{W_0}(w_q\cdot\widetilde{w}_{2(m-q)-1})\\ &\hspace{1cm}+ b_{1,q-1,2(m-q)}\cdot \textsf{Res}_W^{W_0} (w_q\cdot\widetilde{w}_{2(m-q)}), \end{aligned}\end{equation}
which yields :
\begin{equation}
a'_q=\underset{0\leq j \leq 2(m-q)}{\sum} b_{1,q-1,j}\cdot \textsf{Res}_W^{W_0} (w_q\cdot\widetilde{w}_j).
\end{equation}
This ends the induction and the proof of Proposition \ref{rec}. $\blacksquare$
\end{dem}

We eventually get that there exist some coefficients $C_{i,j}\in H^*(k_0,\mathbb Z/2\mathbb Z)$ such that 
\begin{equation*}\textsf{Res}_W^{W_0}(a)=\underset{0\leq i\leq m-1, 0\leq j\leq 2(m-i)}{\sum} C_{i,j}\cdot \textsf{Res}_W^{W_0}(w_i\cdot \widetilde{w}_j) +
a_{m}\end{equation*} where \begin{equation*} a_m=\underset{ m-1\leq i\leq m-1, 0\leq j\leq 2(m-1-i)}{\sum} b_{1,i,j}\cdot w_1^{\mathbb D_4}\cdot w_i^{W'}\cdot\widetilde{w}_j^{W'}= b_{1,m-1,0}\cdot w_1^{\mathbb D_4}\cdot w_{m-1}^{W'}.\end{equation*}

By Proposition \ref{siw0}, $\textsf{Res}_W^{W_0}(w_m)= w_m^{W'}+w_1^{\mathbb D_4}\cdot w_{m-1}^{W'}$ and $w_m^{W'}=0$. Hence, $\textsf{Res}_W^{W_0}(a)$ is a linear combination of the invariants $\textsf{Res}_W^{W_0}(w_i\cdot \widetilde{w}_j)$ of $W_0$, for $0\leq i\leq m$ and $0\leq j\leq 2(m-i)$. Since the restriction to $W_0$ completely determines the invariant $a$, we get that $a$ is a linear combination of the invariants $w_i\cdot \widetilde{w}_j$ of $W$, $0\leq i\leq m$, $0\leq j\leq 2(m-i)$. Therefore, for $0\leq i\leq m, 0\leq j\leq 2(m-i)$, the invariants $w_i\cdot\widetilde{w}_j$  of $W$ generate the module $\textsf{Inv}_{k_0}(W,\mathbb Z/2\mathbb Z)$. This ends the proof of Proposition \ref{gen}.

\subsubsection{A basis of $\textsf{Inv}_{k_0}(W,\mathbb Z/2\mathbb Z)$}
Note that any result of this section is still true if we do not assume that $-1,2\in k_0^{\times 2}$.

\begin{thm}
The family $\{w_i\cdot\widetilde{w}_j\}_{0\leq i\leq m, 0\leq j\leq 2(m-i)}$ is free over $H^*(k_0,\mathbb Z/2\mathbb Z)$.\\
\end{thm}

\begin{dem}
Let $\{\lambda_{i,j}\}_{0\leq i\leq m, 0\leq j\leq 2(m-i)}$ be a family of coefficients of $H^*(k_0,\mathbb Z/2\mathbb Z)$ such that \begin{equation*} a=\underset{0\leq i\leq m, 0\leq j\leq 2(m-i)}{\sum} \lambda_{i,j}\cdot w_i\cdot\widetilde{w}_j=0.\end{equation*}

Let us show by induction on $q\in \{0,...,m\}$, that, for any $q\in\{0,...,m\}$, $\lambda_{q,j}=0$ for any $j\in\{0,...,2(m-i)\}$.\\

Assume first that $q=0$. Let us consider the restriction of $a$ to $H_0$.
We have to show that, for any $0\leq j\leq 2m$, $\lambda_{0,j}=0$.
By Lemma \ref{resh_q}, for any extension $k/k_0$, for any $W$-torsor $T_0$ over $k$ lying in the image of the map $H^1(k,H_0)\rightarrow H^1(k, W)$ and for any $i>0$, we have $w_i(T_0)=0$. Thus, \begin{equation*}a_k(T_0)=\underset{0\leq j\leq 2m}{\sum} \lambda_{0,j}\cdot\widetilde{w}_j(T_0).\end{equation*} 
By Lemma \ref{h0}, the family $\{\textsf{Res}_W^{H_0}(\widetilde{w}_j)\}_{0\leq j\leq 2m}$ is free in the $H^*(k_0,\mathbb Z/2\mathbb Z)$-module of the invariants of $H_0$ modulo $2$. Therefore, we get that $\lambda_{0,j}=0$ for any $0\leq j\leq 2m$.\\

Let now $0<q\leq m$. Let us assume that, for any $0\leq i< q$, $\lambda_{i,j}=0$ for any $0\leq j\leq 2(m-i)$.  
Hence, \begin{equation*}a=\underset{q\leq i\leq m, 0\leq j\leq 2(m-i)}{\sum} \lambda_{i,j}\cdot w_i\cdot\widetilde{w}_j.\end{equation*} 

Let us now consider the restriction $\textsf{Res}_W^{H_q}(a)$ of $a$ to $H_q$. Let $k/k_0$ be an extension and let $T_q=(k(\sqrt {t_1})\times...\times k(\sqrt {t_q})\times k^{2(m-q)}, (u_1,...,u_q,u_{q+1},v_{q+1},...,u_m,v_m))$ be a $W$-torsor over $k$ lying in the image of $H^1(k,H_q)\rightarrow H^1(k,W)$.  By Lemma \ref{resh_q}, if $i\geq q+1$, $w_i(T_q)=0$. Thus, \begin{equation*} a_k(T_q)=\underset{ 0\leq j\leq 2(m-q)}{\sum} \lambda_{q,j}\cdot w_q(T_q)\cdot\widetilde{w}_j(T_q).\end{equation*} Therefore, $\textsf{Res}_W^{H_q}(a)=\underset{ 0\leq j\leq 2(m-q)}{\sum} \lambda_{q,j}\cdot \textsf{Res}_W^{H_q}(w_q\cdot\widetilde{w}_j)$. By Lemma \ref{libHq}, the family $\{\textsf{Res}_W^{H_q}(w_q\cdot\widetilde{w}_j)\}_{0\leq j\leq 2(m-q)}$ is free over $H^*(k_0,\mathbb Z/2\mathbb Z)$. We can conclude that $\lambda_{q,j}=0$ for every $0\leq j\leq 2(m-q)$. This ends the induction. $\blacksquare$
\end{dem}

\subsection{ Cohomological invariants of $W$ of type $B_{2m+1}$}
In this section we will just sketch the proof of Theorem \ref{icbn} with $n$ odd. Let us recall the statement.
\begin{thm*} 
Let $k_0$ be a field of characteristic zero, such that $-1$ and $2$ are squares in $k_0$. Let $m\geq 1$ and let $W$ be a Weyl group of type $B_{2m+1}$. Then the module $\textsf{Inv}_{k_0}(W,\mathbb Z/2\mathbb Z)$ is free over $H^*(k_0,\mathbb Z/2\mathbb  Z)$, with basis \begin{equation*}\{w_i\cdot\widetilde{w}_j\}_{0\leq i\leq m, 0\leq j\leq 2(m-i)}.\end{equation*}
\end{thm*}

Let $W$ be a Weyl group of type $B_{2m+1}$. Let \begin{equation*} S=\{\pm e_i, \pm e_i\pm e_j, 1\leq i\leq 2m+1, 1\leq j\neq i\leq 2m+1 \}\end{equation*} be the root system corresponding to $W$. Let \begin{equation*} S_0=\{\pm e_i, \pm e_i\pm e_j, 1\leq i\leq 2m, 1\leq j\neq i\leq 2m \}\sqcup \{\pm e_{2m+1}\}.\end{equation*} Then the reflections associated with $S_0$ generate a subgroup $W_0$ of $W$, isomorphic to $W'\times \langle r_{e_{2m+1}}\rangle$, where $W'$ is a Weyl group of type $B_{2m}$. 
Let $a\in\textsf{Inv}_{k_0}(W,\mathbb Z/2\mathbb Z)$. Mimicking the proof of the case $B_{2m}$, we may show that $\textsf{Res}_W^{W_0}(a)$ completely determines $a$. Then we write $\textsf{Res}_W^{W_0}(a)$ as cup-products of invariants of $W'$ and of $\mathbb Z/2\mathbb Z$. Looking at the restrictions to the subgroups $H_q$ as in the case $B_{2m}$, we finally get 
the result.  Note that computations are much easier here because of the factor $\mathbb Z/2\mathbb Z$ instead of $\mathbb D_4$.

\section{Cohomological invariants of Weyl groups of type $D_n$}
Note that the case $n=4$ was already treated by Serre in the paper of Knus-Tignol \cite{knus2009}, section 8.\\

Let $n\geq 4$ and let $W$ be a Weyl group of type $D_n$. We associate with $W$ its root system $S=\{\pm e_i\pm e_j\mid 1\leq i<j\leq n\}$ (see Section \ref{weylgroups}). Let us denote by $W'$ the Weyl group of type $B_n$ corresponding to the root system \begin{center}$S'=\{\pm e_i, \pm(e_i\pm e_j) \mid, 1\leq i\leq n, 1\leq j\neq i\leq n\}$.\end{center} We clearly have an inclusion $W\subset W'$. Recall that $W'\simeq \big( \mathbb Z/2\mathbb Z\big)^n\rtimes \mathfrak S_n$. Moreover, $W$ is the kernel of the map \begin{equation*} \begin{aligned}p: \hspace{.5cm} W'\hspace{.5cm} &\rightarrow \mathbb Z/2\mathbb Z\\ \big((\epsilon_1,...,\epsilon_n), \sigma\big)&\mapsto \prod_{i=1}^n \epsilon_i.\end{aligned}\end{equation*}

By Galois descent, we may show that :
\begin{prop}
The image of the map $H^1(k,W)\rightarrow H^1(k, W')$ corresponds to pairs $(L,\alpha)$ such that $\alpha$ has norm $1$ in $L^{\times}/L^{\times 2}$.
\end{prop}

Set $m=[\frac{n}{2}]$. Consider now, as before,  $S_m=\{\pm e_{2i-1}\pm e_{2i}\mid 1\leq i\leq m \}$. It is a root subsystem of $S$. The associated reflections generate the subgroup $H_m$ of $W$. We may show that $H_m$ is a representative of the maximal abelian subgroups of $W$ generated by reflections, up to conjugation (they are indeed all conjugate in $W$). Therefore Theorem \ref{iccox} may be reformulated here as follows :
\begin{cor}
\label{dn}
Let $k_0$ be any field of characteristic zero. The restriction map $\textsf{Res}_W^{H_m}:\textsf{Inv}_{k_0}(W,C)\rightarrow \textsf{Inv}_{k_0}(H_m,\mathbb Z/2\mathbb Z)$ is injective.
\end{cor}

Before to state and prove Theorem \ref{icdn}, let us first describe the cohomological invariants of $H_m$ modulo $2$, fixed by its normalizer in $W$. For any $r,s\geq 0$ and $r+s\leq \frac{n}{2}$ and for any extension $k/k_0$, we set \begin{equation*}\begin{aligned} a_{r,s}: H^1(k,H_m)&\rightarrow H^*(k,\mathbb Z/2\mathbb Z)\\ (u_1,u_1t_1,...,u_n) &\mapsto \underset{\underset{\textrm{ odd numbers}}{1\leq m_1 < ... < m_r \leq n-1 }}{\sum} (u_{m_1})\cdot (u_{m_1+1}) \cdot\ldots\cdot (u_{m_r})\cdot (u_{m_r +1})\\ &\hspace{3cm} \cdot \big( \underset{\underset{  L \cap \{ m_1,m_1+1,...,m_r,m_r+1 \}= \emptyset}{L\in I_s} }{\sum} (u)_L\big)\end{aligned}\end{equation*} 
where \begin{equation*}\begin{aligned} I_s=\{\{l_1,...,l_s\}\in \mathbb N^s \mid 1\leq l_1<...< l_s\leq n \textrm{ and }& \forall 0\leq l\leq \frac{n}{2},\\ &\{2l-1,2l\} \not\subset \{l_1,...,l_s\} \}\end{aligned} \end{equation*} and $(u)_L= (u_{l_1})\cdot\ldots \cdot (u_{l_s})$.
It is easily seen that these maps yield an invariant $a_{r,s}$ of $H_m$. We then get that :
 
\begin{lem}
\label{invH}
The module $\textsf{Inv}_{k_0}(H_m,\mathbb Z/2\mathbb Z)^{N_{H_m}/H_m}$ is free over $H^*(k_0,\mathbb Z/2\mathbb Z)$ with a basis given by invariants $a_{r,s}$, for $ r,s\geq 0$ and $r+s\leq \frac{n}{2}$.
\end{lem}

The proof of Lemma \ref{invH} is left to the reader. Recall now that the Stiefel-Whitney invariants $w_i$ or $\widetilde{w}_i$ of $W'$ give by restriction some Stiefel-Whitney invariants of $W$, still denoted by $w_i$ or $\widetilde{w}_i$ for any $0\leq i\leq n$. We then have Theorem \ref{icdn} :
\begin{thm*}
The $H^*(k_0,\mathbb Z/2\mathbb Z)$-module $\textsf{Inv}_{k_0}(W,\mathbb Z/2\mathbb Z)$ is free with basis $\{w_i\cdot \widetilde{w}_j\}$, where $0\leq i\leq [\frac{n}{2}]$, $0\leq j\leq 2([\frac{n}{2}]-i)$ and $j$ even.
\end{thm*}

\begin{dem}
The details of the computations are left to the reader. For sake of simplicity, let us assume that $-1$ and $2$ are squares in $k_0$. Let us also deal with $n$ even (the case $n$ odd then follows easily) and set $m=\frac{n}{2}$. Let $k/k_0$ be a field extension and let $(L,\alpha)\in H^1(k,W')$ lie in the image of the map $H^1(k,H_m)\rightarrow H^1(k,W')$. Then there exist some $t_1,...,t_m,u_1,...,u_m\in k^\times$ such that $L=k(\sqrt{t_1})\times ... \times k(\sqrt{t_m})$ and $\alpha=(u_1,...,u_m)$. Let $0\leq i\leq m$. We have :
\begin{equation*}
\begin{aligned}
w_i(L,\alpha) =& w_i(\langle 2,2t_1,2,2t_2,...,2,2t_m\rangle)= w_i(\langle t_1,...,t_m\rangle) \\ =& a_{0,i}(u_1,u_1t_1,...,u_m,u_mt_m)
\end{aligned}
\end{equation*}
and
\begin{equation*}
\begin{aligned}
\widetilde{w}_{2j}(L,\alpha) = w_{2j}(\langle u_1,u_1t_1,...,u_m,u_mt_m\rangle)= \sum_{r=0}^{j} a_{r,2(j-r)}(u_1,u_1t_1,...,u_m,u_mt_m)
\end{aligned}
\end{equation*} 
with the convention that $a_{r,s}=0$ if $r+s>m$. We then have, for any $0\leq i\leq m$ and any $0\leq j\leq m-i$ :
\begin{equation*}
\textsf{Res}_W^{H_m}(w_i\cdot\widetilde{w}_{2j})= \sum_{r=0}^j a_{0,i}\cdot a_{r,j-2r}.
 \end{equation*}

A short computation yields that, for $0\leq i\leq m$ and $0\leq j\leq m-i$:
\begin{equation}
\label{eq24}
\textsf{Res}_W^{H_m}(w_i\cdot\widetilde{w}_{2j})= \sum_{r=0}^j  \begin{pmatrix}i+2(j-r) \\ 2(j-r) \end{pmatrix} a_{r,i+2(j-r)}.
\end{equation}

\noindent By using Equation (\ref{eq24}), an easy induction proves that, for any $r,s\geq 0$ such that $r+s\leq m$, $a_{r,s}$ may be written as a linear combination of $\textsf{Res}_W^{H_m}(w_{2r+s-j}\cdot \widetilde{w}_{2j})$, for $0\leq j\leq r$. Therefore, the map $\textsf{Res}_W^{H_m}:\textsf{Inv}_{k_0}(W,\mathbb Z/2\mathbb Z)\rightarrow \textsf{Inv}_{k_0}(H_m,\mathbb Z/2\mathbb Z)^{N_{H_m}/H_m}$ is surjective. $\blacksquare$
\end{dem}

\bibliographystyle{plain}
\bibliography{biblio} 

\end{document}